\documentclass[10pt]{article}
\usepackage{graphicx,color}
\usepackage{amsmath}
\usepackage{amsfonts}
\usepackage{amssymb}

\newtheorem{theorem}{Theorem}

\newtheorem{corollary}{Corollary}

\newtheorem{lemma}{Lemma}

\newtheorem{remark}{Remark}
\newtheorem{conjecture}{Conjecture}
\newenvironment{proof}[1][Proof]{\textbf{#1.} }%
               {\ \rule{0.5em}{0.5em}\medskip}

\newif\ifShowComments
\ShowCommentstrue
\def\strutdepth{\dp\strutbox}
\def\druk#1{\strut\vadjust{\kern-\strutdepth
        {\vtop to \strutdepth{%
                \baselineskip\strutdepth\vss
                        \llap{\hbox{#1}\quad}\null}}}}

\def\bl{\bigl}%
\def\br{\bigr}%

\begin{document}
\title{Entropic repulsion of an interface in an external field}
\author{Y.Velenik\\Laboratoire d'Analyse, Topologie, Probabilit\'es\\UMR-CNRS 6632\\ C.M.I., Universit\'e de Provence\\Marseille, France\\{\tt
velenik@cmi.univ-mrs.fr}}
\maketitle
\begin{abstract}
We consider an interface above an attractive hard wall in the complete wetting regime, and submitted to the action of an external increasing, convex potential, and study its delocalization as the intensity of this potential vanishes. Our main motivation is the analysis of critical prewetting, which corresponds to the choice of a linear external potential.

We also present partial results on critical prewetting in the two dimensional Ising model, as well as a few (weak) results on pathwise estimates for the pure wetting problem for effective interface models.
\end{abstract}

\section{Introduction and results}
%
%
%
%
%
%
\subsection{Effective interface models}
There has been a lot of interest recently in the properties of effective interface models in the presence of an attractive hard wall, and the associated wetting phase transition~\cite{BoDeZe2000,CaVe2000,IsYo2001}, as well as the related problems of entropic repulsion~\cite{BrElFr1986, BoDeZe1995,De1996,BoDeGi2001,DeGi1999,DeGi2000,BeGi2002a, BeGi2002b} and pinning by a local potential~\cite{DuMaRiRo1992,BoBr2001,DeVe2000,IoVe2000,BoVe2001}.

In the present work, we consider $d$-dimensional effective interface models with
convex interactions. That is, the interface is described by a collection of
non-negative real numbers $X_i$, $i\in\Lambda_N=\{-N,\dots,N\}^d$, representing
the height of the interface above site $i$. The probability measure is given by
\begin{multline}
\mathrm{P}_{N,+,\lambda,\upsilon}(\mathrm{d}X) =
\bl(Z_{N,+,\lambda,\upsilon}\br)^{-1}\, \exp\Bigl\{ - \sum_{\langle i,j
\rangle\cap\Lambda_N \neq\emptyset}\!\! \mathsf{U}(X_i-X_j) - \sum_{i\in\Lambda_N}
\lambda\, \mathsf{V}(X_i)\Bigr\}\\
\times \prod_{i\in\Lambda_N} \left( \mathrm{d} X_i + \upsilon \delta_0(\mathrm{d} X_i) \right) \prod_{i\in\mathbb{Z}^d \setminus \Lambda_N} \delta_0(\mathrm{d} X_i)\,,
\label{eq_measure}
\end{multline}
where $\delta_0$ is the Dirac mass at $0$, $\lambda$ is a strictly positive real number, $\upsilon \geq 0$ and the function $\mathsf{U}$ satisfies
$c<\mathsf{U}''<1/c$ for some $c>0$. We are mostly interested in the case $\mathsf{V}(x)=|x|$, but we allow for general $\mathsf{V}:\mathbb{R}^+ \to \mathbb{R}^+$, which are convex, increasing, and satisfy $\mathsf{V}(0)=0$ and the following growth condition:
There exists $f:\mathbb{R}^+\to\mathbb{R}^+$ such that, for all $\alpha>0$,
$$
\limsup_{x\to\infty} \frac {V(\alpha x)} {V(x)} \leq f(\alpha) < \infty\,.
$$
(Obviously any convex polynomial is fine.)
The parameter $\upsilon$ corresponds to the
strength of the pinning potential locally attracting the interface to the wall.
\begin{remark}
Our results and proofs can be extended straightforwardly to the case of a square-well potential of the form $-b\mathbf{1}_{X_i\leq a}$, $b\geq 0$, $a>0$.
\end{remark}
When the pinning potential is absent (i.e.
$\upsilon=0$), we simply omit the corresponding subscript. We denote by
$\upsilon_{\mathrm{c}}$ the critical value of the pinning parameter,
\begin{equation}
\label{eq_CW}
\upsilon_{\mathrm{c}} = \sup\left\{ \upsilon \,:\, \lim_{N\to\infty} |\Lambda_N|^{-1} \log
\frac{Z_{N,+,0,\upsilon}} {Z_{N,+,0}} = 0 \right\}\,.
\end{equation}
This corresponds to the value of the pinning parameter at which the wetting transition takes place. In the case
$\mathsf{U}(x)=x^2$, it is known that $\upsilon_{\mathrm{c}}>0$ when $d\leq 2$~\cite{CaVe2000}, while $\upsilon_{\mathrm{c}} = 0$ when $d\geq
3$~\cite{BoDeZe2000}. (The situation is different when one drops the assumption
of strict convexity of $\mathsf{U}$, see~\cite{CaVe2000}.) We write $\mathsf{CW} = \{\upsilon\geq 0 \,:\, \upsilon < \upsilon_c \} \cup \{0\}$ to denote the region of complete wetting; for these values of the pinning parameter, the interface is delocalized when $\lambda=0$.

Let $\mathcal{H}_{\lambda,d}$ be defined as follows: In dimension $d=2$, $\mathcal{H}(\lambda;2) = |\log\lambda|$; in dimensions $d\geq 3$,  $\mathcal{H}_{\lambda,d} = |\log\lambda|^{1/2}$. Then our main result can be stated as
\begin{theorem}\label{thm_EffInt} 
Let $\upsilon\in\mathsf{CW}$. There exist dimension-dependent, strictly positive constants $\delta$, $\lambda_0$, $r$, and $c$ such that: For all $\lambda<\lambda_0$ and $N> \lambda^{-r}$,
$$
\mathrm{P}_{N,+,\lambda,\upsilon}\bigl(|\Lambda_N|^{-1}\sum_{i\in\Lambda_N} X_i \not\in (\delta \mathcal{H}_{\lambda,d},
\delta^{-1} \mathcal{H}_{\lambda,d})\bigr) \leq \exp\left\{-c\,\lambda\, \mathsf{V}\left(\mathcal{H}_{\lambda,d}\right)\, |\Lambda_N| \right\}\,.
$$
Actually, the upper bound can be strengthened as
$$
\mathrm{P}_{N,+,\lambda,\upsilon}\bigl(|\Lambda_N|^{-1}\sum_{i\in\Lambda_N} X_i \not\in (\delta \mathcal{H}_{\lambda,d},
\delta^{-1} \mathcal{H}_{\lambda,d})\bigr) \leq \exp\left\{-c\,\delta^{-1}\, \lambda\, \mathsf{V}\left(\mathcal{H}_{\lambda,d}\right)\, |\Lambda_N| \right\}\,.
$$
\end{theorem}
\begin{remark}
We do not state and prove here the corresponding results in dimension $1$. The reason is that, in that case, much more detailed informations, valid for an immensely larger class of interactions $\mathsf{U}$ can be obtained. This is of interest, since one would like to understand the degree of universality in the behavior found here. These results will appear elsewhere~\cite{HrVe2003}. Let us just mention that, when $d=1$, the interface is repelled to a height of order $\mathcal{H}$, the unique solution of the equation
$$
\lambda \mathcal{H}^2 \mathsf{V}(2\mathcal{H}) = 1\,.
$$
\end{remark}
The previous theorem only yields estimates on the height of the field averaged
over the box. It is however very easy to derive from it the following local
estimates.
\begin{corollary}
\label{cor_highdim}
Let $\upsilon\in\mathcal{CW}$ and fix $0<\epsilon<1$.
There exist $\delta>0$ and $\lambda_0>0$ such that, for
all $\lambda<\lambda_0$ and $N>N_0(\lambda)$,
\begin{align*}
\sup_{i\in\Lambda_N}\mathrm{E}_{N,+,\lambda,\upsilon}(X_i) &\leq \delta^{-1}|\log\lambda|^{\nu(d)},\\
\inf_{i\in\Lambda_{\epsilon N}}\mathrm{E}_{N,+,\lambda,\upsilon}(X_i) &\geq
\delta\,|\log\lambda|^{\nu(d)}\,;
\end{align*}
\end{corollary}

\bigskip
Our main motivation for studying this model comes from the physical phenomenon of \textit{critical prewetting}.
Consider
some substance in the regime of phase coexistence; let us call the two equilibrium phases
$A$ and $B$. Suppose that phase $A$ occupies the bulk of the vessel, while
the boundary of the latter favors phase $B$. As a result, a film of $B$ phase
is generated along the walls. It is well-known that by increasing the
temperature, the system may undergo a surface phase transition, the wetting
transition, reflected in the behavior of this film: Below the wetting
temperature the film has a microscopic width (partial wetting), while above it
its width becomes macroscopic (complete wetting), in the sense that it diverges
in the thermodynamic limit.

Suppose now that, starting from the complete wetting regime, the system is
pulled away from phase coexistence and only phase $A$ remains thermodynamically
stable. The film of (now unstable) $B$ phase, though still present, cannot
occupy a macroscopic region and therefore its width stays microscopic (i.e.
remains finite in the thermodynamic limit). Critical prewetting corresponds to
the (continuous) divergence of the film thickness as the system nears the phase
coexistence manifold. If we denote by $\lambda$ the free energy difference
between phase $B$ and phase $A$, then one is interested in the behavior of the
average film width as a function of $\lambda$.

The Gibbs measures~\eqref{eq_measure} provide a natural modelization of this phenomenon. Indeed, when the external potential is chosen as $\mathsf{V}(x)=x$, these measures introduce a penalization of the form ``$\lambda \times$ the volume of the wet layer'', which correspond to the free energetic cost for creating such a layer. In that case, the averaged width of the wetting layer is of order $\lambda^{-1/3}$ when $d=1$~\cite{HrVe2003}.

In the Subsection~\ref{ssec_IsingIntro}, we introduce another, more realistic modelization of this phenomenon in a lattice gas, the two-dimensional Ising model, and prove some partial results indicating that this model has the same type of critical behavior, in particular that the critical exponent is also $1/3$.
%
%
\subsection{Some additional results about effective interface models}
In this section, we present some new results about effective interface models, not directly related to the main problem investigated in the paper. Some turn out to be useful in the proof of Theorem~\ref{thm_EffInt}, other are just of independent interest.

\bigskip
The first problem concerns the probability that an interface stays inside a fixed horizontal slab. This problem was first investigated in~\cite{BrElFr1986} in the Gaussian case. Our first result in an extension to the case of strictly convex interactions; our proof is also completely different. Let $\mathrm{P}_{N}$ denote the Gibbs measure with 0-b.c. without the positivity constraint and with $\lambda=\eta=0$. 

\begin{theorem}\label{thm_Slab}
There exist strictly positive constants
$c_-(d)$, $c_+(d)$, $c$ and $\ell_0<\infty$ such that, for any $\ell>\ell_0$ the following holds:
\begin{itemize}
\item If $d=2$, for all $N>e^{c\ell}$,
$$
e^{-e^{-c_+ \ell}\, |\Lambda_N|}
\geq \mathrm{P}_{N}(|X_i| \leq \ell,\, \forall i \in\Lambda_N) 
\geq e^{-e^{-c_- \ell}\, |\Lambda_N|}\,.
$$
\item If $d\geq3$, for all $N>e^{c\ell^2}$,
$$
e^{-e^{-c_+ \ell^2}\, |\Lambda_N|}
\geq \mathrm{P}_{N}(|X_i| \leq \ell,\, \forall i \in\Lambda_N) 
\geq e^{-e^{-c_- \ell^2}\, |\Lambda_N|}\,.
$$
\end{itemize}
\end{theorem}
We have not been able to extend the results on the variance and the mass proved in~\cite{BrElFr1986} in the Gaussian setting. The corresponding results in dimension $1$, valid for a much larger class of models will appear in~\cite{HrVe2003}

\medskip
All other results concern the (pure) wetting transition in effective interface models. We start with results giving more informations on pathwise properties.
\begin{theorem}\label{thm_Wetting}
\begin{enumerate}
\item Let $\upsilon\in\mathsf{CW}$. For all $\epsilon>0$, there exist $N_0(\upsilon,\epsilon)$ and $c(\upsilon)>0$ such that, for all $N\geq N_0$,
$$
\mathrm{P}_{N,+,0,\upsilon} (|\Lambda_N|^{-1}\,\sum_{i\in\Lambda_N} 1_{\{X_i = 0\}} \geq \epsilon) \leq e^{-c\,\epsilon\, |\Lambda_N|}\,.
$$
\item Let $d=2$ and $\upsilon\in\mathsf{CW}$.
$$
\lim_{N\to\infty} \mathrm{E}_{N,+,0,\upsilon} (|\Lambda_N|^{-1}\, \sum_{i\in\Lambda_N}X_i) = \infty\,.
$$
Moreover, for all $\eta$ small enough, there exists $C>0$ such that
$$
\mathrm{E}_{N,+,0,\upsilon} (|\Lambda_N|^{-1}\, \sum_{i\in\Lambda_N}X_i) \geq C\log N\,.
$$
\item We consider the Gaussian case $\mathsf{U}(x) = x^2/2$. Let $d=2$ and $\upsilon$ be large enough. Then there exists $C<\infty$ such that, for all $i\in\Lambda_N$,
$$
\mathrm{E}_{N,+,0,\upsilon} (X_i) \leq C\,.
$$
Moreover, there exists $C>0$ such that, for all $x\in\mathbb{S}^1$,
$$
\lim_{k\to\infty}-\frac 1{k} \limsup_{N\to\infty}\log\mathrm{cov}_{N,+,0,\upsilon} (X_0,X_{[kx]}) > C\,,
$$
where $[y]\in\mathbb{Z}^d$ denotes the componentwise integer part of $y\in\mathbb{R}^d$.
\end{enumerate}
\end{theorem}
Finally, a long-term goal would be to control quantitatively the divergence of the interface as the wetting transition is approached from the partial wetting regime, i.e. the limit $\upsilon\searrow\upsilon_{\rm c}$. This seems far from what can be achieved with today's techniques, but at least we can show (easily), in the $d\geq 3$ Gaussian case, that the transition is second order, in the sense that the average height does indeed diverge as $\upsilon\searrow\upsilon_{\rm c}=0$, and not stay finite and jump at the transition.
\begin{theorem}\label{thm_2ndOrder}
Let $d\geq 3$ and suppose that $\upsilon_{\rm c}=0$ (this is the case, e.g., when $\mathsf{U}(x)=x^2/2$). Then, for any $i\in\mathbb{Z}^d$,
$$
\lim_{\upsilon\searrow\upsilon_{\rm c}} \lim_{N\to\infty} \mathrm{E}_{N,+,\lambda,\upsilon} (X_i) =\infty\,.
$$
\end{theorem}
%
%
\subsection{The two-dimensional Ising model}
\label{ssec_IsingIntro}
It would of course be very interesting to obtain results similar to those of Theorem~\ref{thm_EffInt} also in the case of lattice gases. However, such systems do not always display critical prewetting. In particular, below the roughening temperature, e.g. in the low-temperature $3$-dimensional Ising model, or in subcritical Ising models in dimensions $d\geq 4$, it is expected that the continuous divergence of the width of the wet layer is replaced by an infinite sequence of first order phase transitions, the \textit{layering transitions}, at which the interface jumps up one microscopic unit. This has never been established for the Ising model, though that should follow from a rather straightforward, but technically involved, Pirogov-Sinai analysis. This problem was studied in detail in the series of papers~\cite{DiMa1994,CeMa1996,LeMa1996} in the case of the discrete SOS model.

On the other hand, critical prewetting is expected to occur in the rough phase of such models. In particular, it should occur for any subcritical temperatures in the case of the two-dimensional Ising model. This is what we analyze in this section.

\bigskip
Let $N$ be even and set $\Lambda_N=\{-N/2+1,\dots,N/2\} \times \{0,\dots,N-1\}$,
$\partial^+\Lambda_N = \{x\in\Lambda_N\,:\, \exists y\not\in\Lambda_N,\,
|x-y|=1,\, y_2\geq 0\}$ and $\partial^-\Lambda_N = \{x\in\Lambda_N\,:\, \exists
y\not\in\Lambda_N,\, |x-y|=1,\, y_2=-1\}$.

We define the following Hamiltonian acting on configurations
$\sigma\in\{-1,1\}^{\Lambda_N}$,
$$
H_{\lambda,h,\Lambda_N} (\sigma) = -\sum_{\langle x,y \rangle \subset
\Lambda_N} \sigma_x\sigma_y - \lambda \sum_{x\in\Lambda_N} \sigma_x -
\sum_{x\in\partial^+\Lambda_N} \sigma_x + h \sum_{x\in\partial^-\Lambda_N}
\sigma_x\,.
$$
$\lambda>0$ is the bulk magnetic field, while the boundary field $h>0$ serves to
model the preference of the bottom wall toward $-$ spins. The Gibbs
measure in $\Lambda_N$ is the probability measure on
$\sigma\in\{-1,1\}^{\Lambda_N}$ given by
$$
\mu_{\beta,\lambda,h,\Lambda_N} (\sigma) = (Z_{\beta,\lambda,h,\Lambda_N})^{-1}\;
e^{-\beta\, H_{\lambda,h,\Lambda_N} (\sigma)}\,.
$$

Let $\beta>\beta_c$ and $\lambda=0$. Because of our choice of boundary
conditions, the bulk of the system is occupied, with probability asymptotically
$1$, by the $+$ phase, while the introduction of the boundary field $h$ results
in the creation of a layer of $-$ phase along the bottom wall. It is
well-known~\cite{Ab1980, FrPf1987a} that there is a critical value $h_w(\beta)>0$
of $h$ at which a wetting transition occurs. For $h<h_w(\beta)$, the layer has
a finite thickness (see~\cite{PfVe1997a}), while for $h\geq h_w(\beta)$, the width of the film of $-$
phase becomes macroscopic (of order $N^{1/2}$) (this can be proved, for example, along the line of the argument leading to~\eqref{7_8} below). It is actually expected that in this regime, once suitably rescaled, the microscopic interface converges to Brownian excursion; the sketch proof of such a result in the case $h=1$ can be found in~\cite{Do1993}.

Let us now choose some $h\geq h_w(\beta)$ and let $\lambda>0$.
Our second result describes the behavior of the film of $-$
phase as $\lambda \to 0^+$. There is no canonical way of measuring the width of
the film of $-$ phase in the Ising model, or, for that matter, even to define
the film itself. We proceed as follows.

\begin{figure}[t!]
\centerline{\includegraphics[height=45mm]{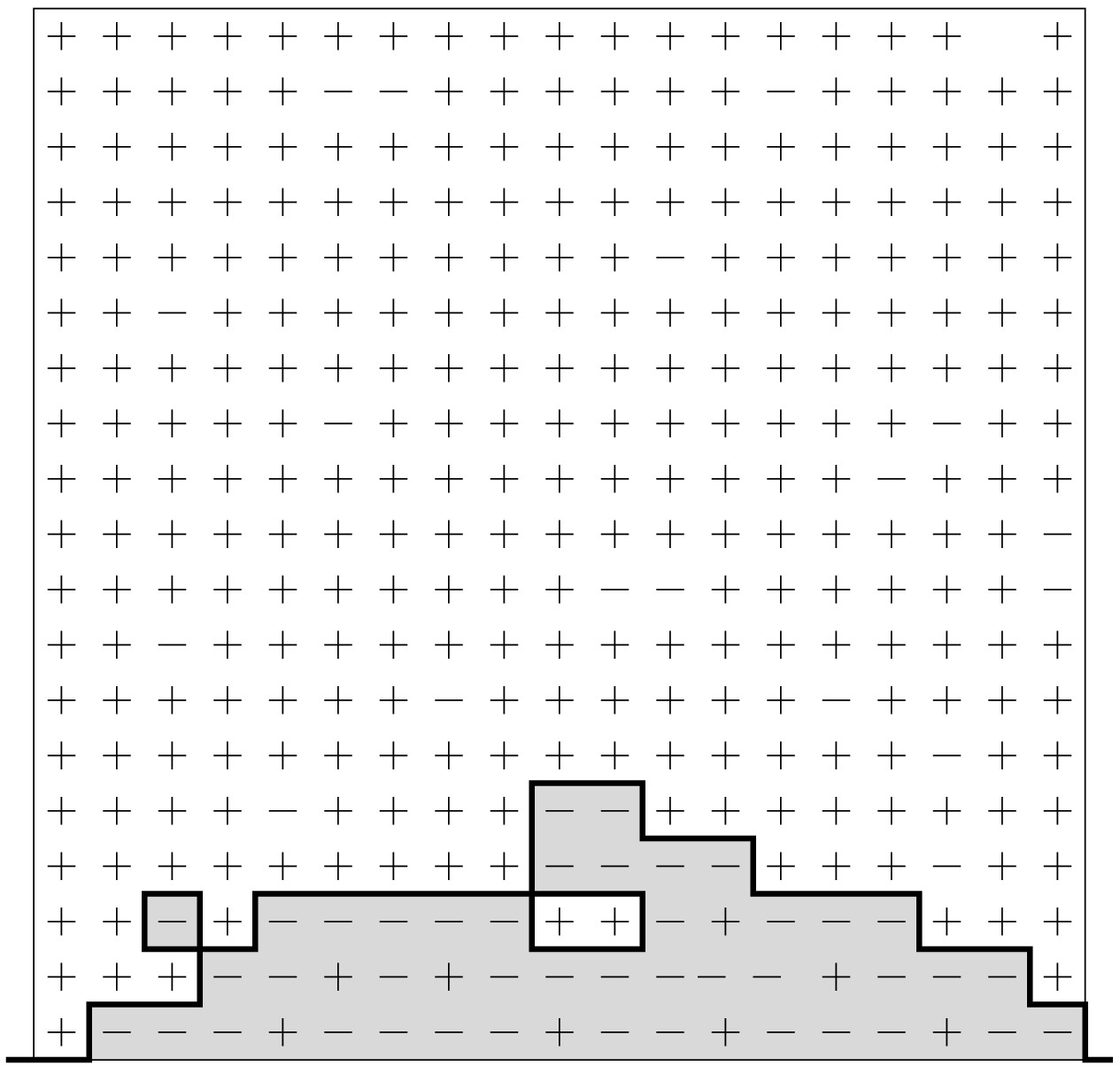}}
\caption{Part of the open contour $\gamma(\sigma)$ (bold line), and the set
$\Lambda^-(\sigma)$ (shaded).}
\label{fig_layervol}
\end{figure}

Given a configuration $\sigma$ in the box $\Lambda_N$, construct its extension
$\bar\sigma$ to $\mathbb{Z}^2$ by setting, for $x\not\in\Lambda_N$,
$\bar\sigma_x = -1$ if $x_2<0$ and $\bar\sigma_x = 1$ otherwise. Then in the
configuration $\bar\sigma$ there is a single (infinite) open Peierls contour.
This Peierls contour can intersect itself, and we use the deformation rules
given in Fig.~\ref{fig_rules} to turn it into a simple curve. We call the
resulting open contour $\gamma(\sigma)$. This contour splits $\mathbb{Z}^2$
into two components. We denote by $\Lambda^-(\sigma)$ the set of all sites of
$\Lambda_N$ lying ``below'' $\gamma(\sigma)$, see Fig.~\ref{fig_layervol}. The
volume of $\Lambda^-(\sigma)$ is one of the measures of the size of the layer
we use. The other one is $|C^-(\sigma)|$, where $C^-(\sigma)$ is the set of all sites in $\Lambda^-(\sigma)$ connected to the bottom wall by a path of $-$ spins in $\sigma$. 

Our result is the following
\begin{theorem}
\label{thm_Ising}
Let $\beta>\beta_c$ and $h\geq h_w(\beta)$. There exist $\lambda_0>0$,
$K<\infty$, and $C_2>0$ such that, for all $0<\lambda<\lambda_0$ and
$N>K\lambda^{-2/3}|\log\lambda|^{3}$,
\begin{equation}\label{eq_Ising_LB}
\mu_{\beta,\lambda,h,\Lambda_N} (|\Lambda^-| <
|\log\lambda|^{-3}\, \lambda^{-1/3}\, N) \leq
e^{-C_2\,|\log\lambda|^2 \lambda^{2/3}\, N}\,,
\end{equation}
and
\begin{equation}\label{eq_Ising_UB}
\mu_{\beta,\lambda,h,\Lambda_N} (|C^-| >
K\, |\log\lambda|^2\, \lambda^{-1/3}\, N) \leq
e^{-C_2\,|\log\lambda|^2 \lambda^{2/3}\, N}\,.
\end{equation}
\end{theorem}
We expect that the two quantities $|\Lambda^-|$ and $|C^-(\sigma)|$ are
typically the same, up to a multiplicative constant:
\begin{conjecture} 
There exists a constant $c>0$ such that
$$
\lim_{N\to\infty} \mu_{\beta,\lambda,h,\Lambda_N} (|C^-| > c
|\Lambda^-|) = 1\,.
$$
\end{conjecture}
(Of course $|\Lambda^-(\sigma)|\geq |C^-(\sigma)|$ for all $\sigma$.) We were
unable to prove this conjecture. The difficulty is that the system can create
big droplets of $+$ phase inside the layer of $-$ phase, and it seems that a
rather complicated metastability analysis is required. Such an analysis is well
understood in the case of regular macroscopic boxes, see e.g.~\cite{ScSh1996}.
However, we are here in a situation where the box is random, and very
anisotropic.

\begin{remark}
The logarithmic corrections present in Theorem~\ref{thm_Ising} are artificial and should be removed. The reason of their appearance is the lack of local CLT type results for Ising random-lines in the vicinity of the boundary of the system. Such estimates have been proved in~\cite{CaIoVe2003a} in the case of random-lines in the bulk but their extension seems delicate. 
\end{remark}

\vskip 0.2cm
\noindent
{\bf Acknowledgments.} It is a pleasure to thank Ostap Hryniv for many fruitful discussions, as well as useful comments and criticisms on earlier versions of this text. The author is also grateful to Walter Selke for bringing his attention to the problem of interfacial adsorption, the understanding of which was the original motivation for this work.

\section{Proofs}
%
%
%
%
%
%
\subsection{Proof of Theorem~\ref{thm_EffInt} and its corollary}
\subsubsection{Proof of Theorem~\ref{thm_EffInt}}
Let $\nu(2)=1$ and $\nu(d)=1/2$ when $d\geq 3$.
Theorem~\ref{thm_EffInt} is a rather simple corollary of Theorem~\ref{thm_Slab} and the following result.
\begin{lemma}[\cite{DeGi2000}, Lemma~2.9] There exists $K>0$ such that, for all $\epsilon>0$,
\begin{equation}
\label{eq_entrep}
\lim_{N\to\infty} \mathrm{P}_{N,+,0} \left(\left\vert \left\{
i \in \Lambda_{[N/2]} \,:\, X_i \leq K(\log
N)^{\nu(d)} \right\}\right\vert \geq \epsilon |\Lambda_{[N/2]}| \right) = 0\,.
\end{equation}
\end{lemma}

\bigskip
We first prove a lower bound on the partition function.
\begin{lemma} \label{lem_lowerboundPF_highdim}
There exist $\lambda_0>0$, $r>0$ and $C>0$ such that, for all $\lambda<\lambda_0$ and
$N>\lambda^{-r}$, and any nonnegative $\upsilon$,
$$
Z_{N,+,\lambda,\upsilon} \geq e^{-C\, \lambda
\mathsf{V}(|\log\lambda|^{\nu(d)})\, |\Lambda_N|}\, Z_{N,+,0,\upsilon}\,.
$$
\end{lemma}
\begin{proof}[Proof of Lemma~\ref{lem_lowerboundPF_highdim}]
First observe that it is enough to prove the corresponding claim without
pinning. Indeed, by FKG,
\begin{align*}
\frac{Z_{N,+,\lambda,\upsilon}} {Z_{N,+,0,\upsilon}}
&=
\langle e^{-\lambda \sum_i \mathsf{V}(X_i)}
\rangle_{N,+,\upsilon}\\
&\geq
\langle e^{-\lambda \sum_i \mathsf{V}(X_i)}
\rangle_{N,+,0}
=
\frac{Z_{N,+,\lambda}}{Z_{N,+,0}}\,.
\end{align*}
Let now 
$$
\mathcal{H}=\left|\log\left( \lambda \mathsf{V}(|\log\lambda|^{\nu(d)})\right) / c_-(d)
\right|^{\nu(d)}\,,
$$
where $c_-(d)$ has been introduced in Theorem~\ref{thm_Slab} and $\nu(d)$ is defined above.
$$
Z_{N,+,\lambda} \geq e^{-\lambda |\Lambda_N| \mathsf{V}(2\mathcal{H})}\; Z_{N,0}\;
\mathrm{P}_{N,0}(0\leq X_i \leq 2\mathcal{H},\, \forall i \in\Lambda_N)\,,
$$
Obviously, shifting the boundary condition from $0$ to $\mathcal{H}$ and changing
variables, we can write, thanks to the convexity of $\mathsf{U}$,
$$
\mathrm{P}_{N,0}(0\leq X_i \leq 2\mathcal{H},\, \forall i \in\Lambda_N) \geq e^{-C
\mathcal{H}^2\, N^{d-1}} \mathrm{P}_{N,0}(|X_i| \leq \mathcal{H},\, \forall i
\in\Lambda_N)\,.
$$
Applying~Theorem~\ref{thm_Slab}, we get that
$$
Z_{N,+,\lambda} \geq e^{-(\lambda \mathsf{V}(2\mathcal{H}) + e^{-c_- \mathcal{H}^{-\nu(d)}} +
c\mathcal{H}^2N^{-1}) |\Lambda_N|}\; Z_{N,0}\,.
$$
Making use of the growth condition on $\mathsf{V}$, the claim follows.
\end{proof}

Using the convexity of $\mathsf{U}$, the lemma immediately implies that
\begin{align*}
\mathrm{P}_{N,+,\lambda,\upsilon} &\left(\sum_{i\in\Lambda_N} X_i \geq \delta^{-1}
|\log\lambda|^{\nu(d)}\, |\Lambda_N|\right) \\
&\leq
\mathrm{P}_{N,+,\lambda,\upsilon} \left(\sum_{i\in\Lambda_N} \mathsf{V}(X_i) \geq \mathsf{V}\left(\delta^{-1}
|\log\lambda|^{\nu(d)}\right)\, |\Lambda_N|\right)\\
&\leq
\exp\left\{-\lambda \mathsf{V}\left(\delta^{-1}
|\log\lambda|^{\nu(d)}\right)\, |\Lambda_N|  \right\} \; \frac{Z_{N,+,0}}{Z_{N,+,\lambda}}\\
&\leq \exp\left\{-\lambda \left(\mathsf{V}\left(\delta^{-1}
|\log\lambda|^{\nu(d)}\right) - C\, \mathsf{V}\left(|\log\lambda|^{\nu(d)}\right) \right)\, |\Lambda_N|  \right\}\\
&\leq \exp\left\{-\lambda \mathsf{V}(|\log\lambda|^{\nu(d)}) \left(\delta^{-1} - C)  \right)\, |\Lambda_N|  \right\}\,,
\end{align*}
and one half of Theorem~\ref{thm_EffInt} is proved.

Let us consider the second half. We first apply
Lemma~\ref{lem_lowerboundPF_highdim} to remove the external field,
\begin{multline*}
\mathrm{P}_{N,+,\lambda,\upsilon} (\sum_{i\in\Lambda_N} X_i \leq \delta
|\log\lambda|^{\nu(d)})\\
\leq
e^{C\lambda|\log\lambda|^{\nu(d)} N^d}\; \mathrm{P}_{N,+,0,\upsilon}
(\sum_{i\in\Lambda_N} X_i \leq \delta |\log\lambda|^{\nu(d)})\,.
\end{multline*}
Let $R=R_0\cdot\left(\lambda \mathsf{V} (|\log\lambda|^{\nu(d)}) \right)^{-1/d}$, with $R_0>0$ some small
enough real number to be chosen below. Let also
$$
\Delta_R = \{x\in\mathbb{Z}^d \,:\, x_i \equiv 0 \mod 2(R+1),\text{ for some
}i\in\{1,\dots,d\} \}\,.
$$
Using the usual expansion on the pinned sites (see e.g.~\cite{BoVe2001}) and FKG inequalities, we can write
\begin{align*}
\mathrm{P}_{N,+,0,\upsilon}& \bigl(\sum_{i\in\Lambda_N} X_i \leq
\delta |\log\lambda|^{\nu(d)} |\Lambda_N|\bigr)\\
&=
\sum_{A\subset\Lambda_N} \nu_{N,+,\upsilon}(A) \, \mathrm{P}_{A^c,+,0}
\bigl(\sum_{i\in\Lambda_N} X_i \leq \delta |\log\lambda|^{\nu(d)} |\Lambda_N| \bigr)\\
&\leq
\sum_{A\subset\Lambda_N} \nu_{N,+,\upsilon}(A) \, \mathrm{P}_{A^c,+,0} \bigl(\sum_{i\in\Lambda_N}
X_i \leq \delta |\log\lambda|^{\nu(d)} |\Lambda_N|\,\bigm\vert\, X\equiv 0 \text{ on
}\Delta_R \bigr)\,, 
\end{align*}
where $\nu_{N,+,\upsilon}$ is the induced probability measure on the subsets of $\Lambda_N$.

The grid $\Delta_R$ splits $\Lambda_N$ into $M=|\Lambda_N|/(2R+2)^d$ cubic cells of sidelength $2R+1$.
(To simplify the notation, we suppose that this splitting can be done exactly.)
Given a realization $A$ of the random set
of the pinned sites, we declare a cell \textsl{clean} if it does not contain
any site of $A$; otherwise it is \textsl{dirty}.
We need to ensure that, for all $\upsilon\in\mathsf{CW}$, at least half of the cells of the grid are clean. 
This follows easily from the first claim of Theorem~\ref{thm_Wetting} with $\epsilon = \tfrac12 M/|\Lambda_N|$, which implies that at least one half of the cells are clean once $N$ is large enough, up to an event of probability $e^{-c\, M}$.
In the rest of this proof, we suppose $N$ is at least that large.

Let us number the clean cells as $\mathcal{C}_1,\dots,\mathcal{C}_K$,
$K\geq\tfrac12 M$, and let $Y_k$, $k=1,\dots,K$, be the indicator
function of the event
$$
\sum_{i\in \mathcal{C}_k} X_i \leq 4\delta |\log\lambda|^{\nu(d)}\;
|\mathcal{C}_k|\,.
$$
We clearly have the inclusion (for all $\lambda$ small enough)
$$
\bigl\{ \sum_{i\in\Lambda_N} X_i \leq \delta |\log\lambda|^{\nu(d)} |\Lambda_N| \bigr\} \subset
\bigl\{ \sum_{k=1}^K Y_i \geq \tfrac13 K \bigr\}\,.
$$
Now, under $\mathrm{P}_{A^c,+,0} (\,\cdot\,\vert\, X\equiv 0 \text{ on
}\Delta_R)$, the $Y_i$ are i.i.d. Bernoulli random variables. Let us denote
their common law by $\mathrm{Q}_R$. Using~\eqref{eq_entrep}, we see that
\begin{align*}
\mathrm{Q}_R (Y_i = 1)
&\leq
\mathrm{P}_{R,+,0} \left(\left\vert \left\{
i \in \Lambda_{[R/2]} \,:\, X_i \leq 2^{d+3}\delta|\log
\lambda|^{\nu(d)} \right\}\right\vert \geq \tfrac12 |\Lambda_{[R/2]}| \right)\\
&\leq 1/10\,,
\end{align*}
provided $\delta$, and then $\lambda$, are chosen small enough. Consequently,
$$
\otimes_{k=1}^M \mathrm{Q}_R (\sum_{k=1}^M Y_i \geq \tfrac16 M) \leq e^{-C\,
M}\,,
$$
and the conclusion follows, provided we take $R_0$ small enough.

\subsubsection{Proof of Corollary~\ref{cor_highdim}}
We only give the proof for $d=2$ since the argument for $d\geq 3$ is identical.
Notice first that it is enough to prove the result for $X_0$. Indeed, by FKG,
$$
\mathrm{E}_{N,+,\lambda,\upsilon}(X_i) \leq \mathrm{E}_{2N,+,\lambda,\upsilon}(X_0)\,,
$$
for all $i\in\Lambda_N$, while for $i\in\Lambda_{\epsilon N}$
($0<\epsilon<1$) we have that
$$
\mathrm{E}_{N,+,\lambda,\upsilon}(X_i) \geq \mathrm{E}_{\epsilon N,+,\lambda,\upsilon}(X_0)\,,
$$
Theorem~\ref{thm_EffInt} implies that
$$
\mathrm{E}_{N,+,\lambda,\upsilon}(|\Lambda_N|^{-1} \sum_{i\in\Lambda_N} X_i ) \in
(a_1,a_2) \mathcal{H}_{\lambda,d}\,,
$$
for some constants $0<a_1<a_2<\infty$ provided $\lambda$ is small enough and $N$
large enough. Using FKG inequalities, we deduce that
$$
|\Lambda_N|^{-1} \sum_{i\in\Lambda_N} \mathrm{E}_{N,+,\lambda,\upsilon}(X_i)
\leq
|\Lambda_N|^{-1} \sum_{i\in\Lambda_N} \mathrm{E}_{2N,+,\lambda,\upsilon}(X_0)
=
\mathrm{E}_{2N,+,\lambda,\upsilon}(X_0)\,,
$$
from which it follows that $\mathrm{E}_{N,+,\lambda,\upsilon}(X_0) \geq a_1
\mathcal{H}_{\lambda,d}$ for $N$ large enough. The complementary bound is also a trivial
consequence of FKG:
\begin{align*}
|\Lambda_N|^{-1} \sum_{i\in\Lambda_N} \mathrm{E}_{N,+,\lambda,\upsilon}(X_i)
&\geq
|\Lambda_N|^{-1} \sum_{i\in\Lambda_{[N/2]}} \mathrm{E}_{N,+,\lambda,\upsilon}(X_i)\\
&\geq
|\Lambda_N|^{-1} \sum_{i\in\Lambda_{[N/2]}} \mathrm{E}_{[N/2],+,\lambda,\upsilon}(X_0)\\
&=
\tfrac14\; \mathrm{E}_{[N/2],+,\lambda,\upsilon}(X_0)\,.
\end{align*}
This shows that $\mathrm{E}_{N,+,\lambda,\upsilon}(X_0) \leq 4 a_2\,  \mathcal{H}_{\lambda,d}$.

%
%
%
%
%
%
\subsection{Proof of Theorems~\ref{thm_Slab}, \ref{thm_Wetting} and \ref{thm_2ndOrder}}
\begin{proof}[Proof of Theorem~\ref{thm_Slab}]
We first establish the lower bound. We write
$$
\mathrm{P}_{N}(|X_i| \leq \ell,\, \forall i \in\Lambda_N) \geq
\mathrm{P}_{N,+}(X_i \leq \ell,\, \forall i \in\Lambda_N)\,
\mathrm{P}_{N}(X_i \geq 0,\, \forall i \in\Lambda_N)\,.
$$
By~\cite[Theorem~3.1]{DeGi2000}, $\mathrm{P}_{N}(X_i \geq 0,\, \forall i
\in\Lambda_N) \geq e^{-c\, N^{d-1}}$ for some $c>0$ if $d\geq 2$, and can thus be neglected. Let us consider the other one. We introduce $R=e^{\kappa\,\ell}$ if $d=2$, and $R=e^{\kappa\,\ell^2}$ if $d\geq 3$, where $\kappa$ is a small number to be chosen later, and the grid
$$
\Delta_R =
\{x\in\mathbb{Z}^d \,:\, x_i \equiv 0 \mod 2(R+1), \text{ for some
}i\in\{1,\dots,d\}\}\,.
$$
This grid splits
$\Lambda_N$ into $M=|\Lambda_N|/(2R+2)^d$ cubic cells of sidelength $2R+1$. (To simplify the notation, we suppose that this splitting can be done exactly.) Denoting by $\mathrm{P}_{N,+}^t$ the measure with $t$ boundary conditions, FKG inequalities then imply that
\begin{align*}
\mathrm{P}_{N,+}(X_i \leq \ell,\, \forall i \in\Lambda_N)
&\geq \mathrm{P}_{N,+}^\ell(X_i \leq \ell,\, \forall i \in\Lambda_N \,|\, X
\equiv \ell \text{ on } \Delta_R)\\
&= \left( \mathrm{P}_{R}(X_i \geq 0,\, \forall i
\in\Lambda_R \,|\, X_i\leq \ell,\, \forall i\in \Lambda_R)
\right)^{M}\,.
\end{align*}
Now,
\begin{multline*}
\mathrm{P}_{R}(X_i \geq 0,\, \forall i
\in\Lambda_R \,|\, X_i\leq \ell,\, \forall i\in \Lambda_R)\\
\geq\mathrm{P}_{R,+}(X_i\leq \ell,\, \forall i\in \Lambda_R)
\,  \mathrm{P}_{R}(X_i \geq 0,\, \forall i
\in\Lambda_R)\,.
\end{multline*}
As above,
$$
\mathrm{P}_{R}(X_i \geq 0,\, \forall i
\in\Lambda_R) \geq e^{-c\, R^{d-1}}\,,
$$
while, using~\cite[Corollary~2.6]{DeGi2000},
$$
\mathrm{P}_{R,+}(X_i\leq \ell,\, \forall i\in \Lambda_R)
\geq 1 - |\Lambda_R|\,
\sup_{i\in\Lambda_R}\mathrm{P}_{R,+}(X_i > \ell) \geq
\tfrac12\,,
$$
provided that $\kappa$ is chosen small enough, and then $\ell$ large enough
(actually the above reference does not give explicitly a rate of convergence,
but a glance at the proof shows that $\sup_{i\in\Lambda_R}\mathrm{P}_{R,+}(X_i > \ell) \leq R^{-c\,\ell}$).

Collecting all the estimates proves the lower bound. 

\medskip
Let us now turn to the upper bound. $\Delta_R$ is the same grid as above (for an
$R$ to be chosen later). Then it follows from convexity of $\mathsf{U}$ and FKG inequalities that
\begin{align*}
\mathrm{P}_{N}(|X_i| \leq \ell,\, \forall i \in\Lambda_N)
&\leq
\mathrm{P}_{N,+}(X_i \leq 2\ell,\, \forall i \in\Lambda_N)\,
e^{c\ell^2\,|\partial \Lambda_N|}\\
&\leq
\mathrm{P}_{N,+}(X_i \leq 2\ell,\, \forall i \in\Lambda_N \,|\, X\equiv 0
\text{ on $\Delta_R$ })\, e^{c\ell^2\,|\partial \Lambda_N|}\\
&=
\left( \mathrm{P}_{R,+}(X_i \leq 2\ell,\, \forall i \in\Lambda_R)
\right)^{M}\, e^{c\ell^2\,|\partial \Lambda_N|}\,.
\end{align*}
We choose $R$ as follows: If $d=2$, then $R=e^{\rho \ell}$, while in the case $d\geq 3$, we choose $R=e^{\rho \ell^2}$, with $\rho$ large enough.
Then
$$
\mathrm{P}_{R,+}(X_i \leq 2\ell,\, \forall i \in\Lambda_R) \leq 1/2\,,
$$
by the results of~\cite{DeGi2000} (since the typical height of the repelled field is
proportional to  $(\log R)^{\nu(d)} = \rho^{\nu(d)} \ell$). This completes the
proof.
\end{proof}

\noindent
\begin{proof}[Proof of Theorem~\ref{thm_Wetting}]

\noindent
\textit{First statement.}
Let $\mathcal{A}$ be the (random) set of pinned sites.
The crucial observation is that
\begin{align*}
\mathrm{E}_{N,+,0,\upsilon} \left( e^{t|\mathcal{A}|} \right) 
&= \sum_{A\subset\Lambda_N} e^{t|A|}\, \upsilon^{|A|} \mathrm{Z}_{A^c,+,0,0} / \mathrm{Z}_{N,+,0,\upsilon}\\
&= \frac{\mathrm{Z}_{N,+,0,e^t\upsilon}} {\mathrm{Z}_{N,+,0,\upsilon}}\,.
\end{align*}
Therefore an application of the exponential Chebyshev inequality yields
$$
\mathrm{P}_{N,+,0,\upsilon}\left( |\mathcal A| > \zeta |\Lambda_N|  \right) 
\leq e^{-t\zeta|\Lambda_N|  + \log(\mathrm{Z}_{N,+,0,e^t\upsilon}/\mathrm{Z}_{N,+,0}) - \log(\mathrm{Z}_{N,+,0,\upsilon}/\mathrm{Z}_{N,+,0} )}\,.
$$
Choosing $t$ such that $e^t=\tfrac12(\upsilon+\upsilon_c)$ and using~\eqref{eq_CW}, we get the conclusion, for all $N>N_0(\upsilon,\zeta)$.

\noindent
\textit{Second statement.}
Since $\lambda=0$ in all the proof, we omit it from the notations.
It is proved in~\cite{CaVe2000} that the condition $\upsilon\in\mathsf{CW}$ implies that
$$
\lim_{N\to\infty} |\Lambda_N|^{-1}\; \mathbb{E}_{N,+,\upsilon} \left[
\rho_N\right] =0\,,
$$
where $\rho_N(X)=\sum_{i\in\Lambda_N} 1_{\{X_i=0\}}$. So there exists a sequence $c_N$, with $\lim c_N=0$, such that $\mathbb{E}_{N,+,\upsilon}
\left[ \rho_N \right] \leq c_N\, |\Lambda_N|$. Let
$\mathcal{E}_N = \{X\,:\,\rho_N(X)\leq 2c_N\, |\Lambda_N|\}$. Using again the expansion on pinned sites, we have the following 
\begin{align*}
\mathbb{E}_{N,+,\upsilon} \left[
\sum_{i\in\Lambda_N} X_i \right]
&\geq \mathbb{E}_{N,+,\upsilon} \left[
\sum_{i\in\Lambda_N} X_i 1_{\mathcal{E}_N} \right]\\
&= \sum_{\substack{A\subset\Lambda_N\\|A|\leq 2c_N |\Lambda_N|}}
\nu_{N,+,\upsilon}(A)\; \mathbb{E}_{A^c,+,0} (\sum_{i\in\Lambda_N} X_i)\\
&\geq \sum_{\substack{A\subset\Lambda_N\\|A|\leq 2c_N |\Lambda_N|}}
\nu_{N,+,\upsilon}(A)\;
 \sum_{\substack{i\in\Lambda_N\\d(x,A)\geq
c_N^{-1/4}}}\mathbb{E}_{B_{c_N^{-1/4}}(x),+,0} (X_i)\\
&\geq C\, |\log c_N|\, |\Lambda_N|\,,
\end{align*}
where $B_r(x) = \{y \,:\,|y-x|\leq r\}$, and we used positivity of the field and FKG inequalities; the last inequality follows from Markov inequality and the usual entropic repulsion estimate for a square box of radius
$c_N^{-1/4}$ with $0$ b.c., see~\cite{DeGi2000}.

\bigskip
Since, when $\upsilon$ is sufficiently small, we can even choose $c_N\sim N^{-1}$~\cite{CaVe2000}, in that case one has the stronger result
$$
\mathbb{E}^{+,\varepsilon}_N \left[
\sum_{x\in\Lambda_N} \phi_x \right] \geq C\, \log N\, |\Lambda_N|\,.
$$

\noindent
\textit{Third statement.}
Since $\lambda=0$ in all the proof, we omit it from the notations.
Our proof is based on the following estimate: There exist $\eta_0$ and $c>0$ such that, for all $B\subset\Lambda$ and all $\upsilon>\upsilon_0$,
\begin{equation}
\label{eq_PercLargeUps}
\nu_{N,+,\upsilon} (\mathcal{A} \cap B = \emptyset) \leq e^{-c|\log\upsilon| |B|}\,.
\end{equation}
Let us first assume that~\eqref{eq_PercLargeUps} holds and prove the two statements of the Theorem. The first statement is very easy: By the usual decomposition over pinned sites, and FKG inequalities,
\begin{align*}
\mathrm{E}_{N,+,0,\upsilon} (X_i) &= \sum_{R\geq 1} \sum_{\substack{A\subset\Lambda_N\\ d(A,x)=R}} \nu_{N,+,\upsilon}(A) \, \mathrm{E}_{A^{\rm c},+,0,0} (X_i)\\
&\leq \sum_{R\geq 1} \sum_{\substack{A\subset\Lambda_N\\ d(A,x)=R}} \nu_{N,+,\upsilon}(A) \, \sup_{|j-i|=R}\mathrm{E}_{\mathbb{Z}^2\setminus\{j\},+,0,0} (X_i)\,.
\end{align*}
Since $\mathrm{E}_{\mathbb{Z}^2\setminus\{j\},+,0,0} (X_i) \leq C\log R$ uniformly in $j$ such that $|j-i|=R$~\cite{DuMaRiRo1992}, the first bound is proved.

Let us now establish the positivity of the mass. Decomposing again according to the pinned sites, we get
\begin{align*}
\mathrm{E}_{N,+,0,\upsilon} (X_i X_j) 
&= \sum_{\substack{A\subset\Lambda_N \\ i \not\stackrel{A^{\rm c}}{\leftrightarrow}j}} \nu_{N,+,\upsilon}(A) \, \mathrm{E}_{A^{\rm c},+,0,0} (X_i X_j)\,,\\
&\hspace{1cm}+ \sum_{\substack{A\subset\Lambda_N \\ i \stackrel{A^{\rm c}}{\leftrightarrow}j}} \nu_{N,+,\upsilon}(A) \, \mathrm{E}_{A^{\rm c},+,0,0} (X_i X_j)\,,
\end{align*}
where $i \stackrel{A^{\rm c}}{\leftrightarrow}j$ means that $i$ and $j$ are connected by a path containing only sites of $\Lambda_N\setminus A$, and $i\not\stackrel{A^{\rm c}}{\leftrightarrow}j$ is the complementary event.

Suppose $i\not\stackrel{A^{\rm c}}{\leftrightarrow}j$ occurs. Let $\gamma_i$ be the innermost closed path of pinned sites surrounding $j$ (possibly using sites outside from $\Lambda_N$), and let $\gamma_j$ be the analogous object for $j$. If $\gamma_i$ does not surround $j$, then let $\gamma=\gamma_i$, otherwise let $\gamma=\gamma_j$. Observe that under $i\not\stackrel{A^{\rm c}}{\leftrightarrow}j$, $\gamma$ surrounds only one of the two sites, and is the innermost path with this property. Suppose, for definiteness that it surrounds $i$. We then bound the first term as follows, writing $\mathrm{int}(\gamma)$ for the set of site surrounded by $\gamma$ or along $\gamma$,
\begin{align*}
\mathrm{E}_{N,+,0,\upsilon} (X_i X_j \,|\, \gamma)
&\leq \mathrm{E}_{\Lambda_N\setminus\mathrm{int}(\gamma),+,0,\upsilon} (X_j)
\; \mathrm{E}_{\mathrm{int}(\gamma)\setminus\gamma,+,0,0} (X_i)\\
&\leq \mathrm{E}_{\Lambda_N,+,0,\upsilon} (X_j)
\; \mathrm{E}_{\mathrm{int}(\gamma)\setminus\gamma,+,0,0} (X_i)\,.
\end{align*}
Since averaging the last expectation over $\gamma$ yields $\mathrm{E}_{N,+,0,\upsilon} (X_i,\, i\not\stackrel{A^{\rm c}}{\leftrightarrow}j,\, \gamma=\gamma_i)$, we see that
$$
\sum_{\substack{A\subset\Lambda_N \\ i \not\stackrel{A^{\rm c}}{\leftrightarrow}j}} \nu_{N,+,\upsilon}(A) \, \mathrm{E}_{A^{\rm c},+,0,0} (X_i X_j)
- \mathrm{E}_{N,+,0,\upsilon} (X_i) \mathrm{E}_{N,+,0,\upsilon} (X_j) \leq 0\,.
$$
Therefore, it only remains to prove that
$$
\sum_{\substack{A\subset\Lambda_N \\ i \stackrel{A^{\rm c}}{\leftrightarrow}j}} \nu_{N,+,\upsilon}(A) \; \mathrm{E}_{A^{\rm c},+,0,0} (X_i X_j)
$$
is exponentially small in $|i-j|$. This follows from~\eqref{eq_PercLargeUps}. Indeed, this expression is bounded above by
$$
\sum_{R_1\geq 0, R_2\geq 0} \sum_{\substack{A\subset\Lambda_N,\, i \stackrel{A^{\rm c}}{\leftrightarrow}j \\ d(A,x)=R_1, d(A,y)=R_2}} \nu_{N,+,\upsilon}(A)\; \mathrm{E}_{A^{\rm c},+,0,0} (X_i X_j)\,.
$$
Now, proceeding as above,
$$
\mathrm{E}_{A^{\rm c},+,0,0} (X_i X_j) \leq \Bigl(\mathrm{E}_{A^{\rm c},+,0,0} (X_i^2) \mathrm{E}_{A^{\rm c},+,0,0} (X_j^2) \Bigr)^{1/2} \leq C \log R_1 \log R_2\,.
$$
Moreover, it follows from~\eqref{eq_PercLargeUps} that
$$
\sum_{\substack{A\subset\Lambda_N,\, i \stackrel{A^{\rm c}}{\leftrightarrow}j \\ d(A,x)=R_1, d(A,y)=R_2}} \nu_{N,+,\upsilon}(A) \leq e^{-C|\log\upsilon| \;|i-j|} e^{-C|\log\upsilon| \; (R_1{}^2+R_2{}^2)}\,,
$$
provided $\upsilon$ is chosen large enough.
The statement follows.

\medskip
Let us turn now to the proof of~\eqref{eq_PercLargeUps}.
The control of the distribution of pinned sites is done similarly as in~\cite{DeVe2000,IoVe2000,BoVe2001}. Let $B\subset\Lambda_N$. We set $B_0=B$, and define
$$
B_{k+1} = B_k \cup \{x\in\mathbb{Z}^2\,:\, d(x,B_k)=1\}\,.
$$
We then have
$$
\nu_{N,+,\upsilon} (\mathcal{A} \cap B = \emptyset) \leq \sum_{k\geq 0} \nu_{N,+,\upsilon} (\mathcal{A} \cap B_k = \emptyset \,|\, \bar{\mathcal{A}} \cap B_{k+1} \neq \emptyset)\,,
$$
where $\bar{\mathcal{A}} = A\cup(\Lambda_N)^{\rm c}$. It is therefore enough to prove that
$$
\nu_{N,+,\upsilon} (\mathcal{A} \cap B_k = \emptyset \,|\, \bar{\mathcal{A}} \cap B_{k+1} \neq \emptyset) \leq e^{-c |\log\upsilon| \, |B|}\,.
$$
We then write, as in~\cite{DeVe2000,IoVe2000,BoVe2001},
\begin{align*}
\nu_{N,+,\upsilon} (\mathcal{A} \cap B_k = \emptyset \,|\, \bar{\mathcal{A}} \cap B_{k+1} \neq \emptyset) 
&\leq \left\{\sum_{C\subset B_k} \upsilon^{|C|} \inf_{\substack {A\cap B_k=\emptyset \\ \bar{A} \cap B_{k+1} \neq \emptyset}} \frac{Z_{(A\cup C)^{\rm c},+,0}}{Z_{A^{\rm c},+,0}}  \right\}^{-1}\\
&\leq \left\{\upsilon^{|B_k|} \inf_{\substack {A\cap B_k=\emptyset \\ \bar{A} \cap B_{k+1} \neq \emptyset}} \frac{Z_{(A\cup B_k)^{\rm c},+,0}}{Z_{A^{\rm c},+,0}}  \right\}^{-1}\,.
\end{align*}
Let us number the sites of $B_k=\{t_1,\ldots,t_{|B_k|}\}$ in such a way that $t_1$ is a neighbor of a site of $\bar{A}\cap B_{k+1}$, and each $t_k$, $k>1$, is a neighbor of at least one $t_j$, $j<k$. We write the last fraction as a telescoping product: 
$$
\frac{Z_{(A\cup B_k)^{\rm c},+,0}}{Z_{A^{\rm c},+,0}} = \prod_{k=1}^{|B_k|} \frac {Z_{(A\cup \{t_1,\ldots,t_k\})^{\rm c},+,0}} {Z_{(A\cup \{t_1,\ldots,t_{k-1}\})^{\rm c},+,0}}\,.
$$
The conclusion follows once we prove that there exists $c>0$ such that
$$
\frac{Z_{(C\cup \{i\})^{\rm c},+,0}} {Z_{C^{\rm c},+,0}} \geq c\,,
$$
uniformly in $C\subset\Lambda_N$, and $i$ neighboring $C$.
\begin{equation*}
\frac{Z_{(C\cup \{i\})^{\rm c},+,0}}{Z_{C^{\rm c},+,0}} =
\lim_{\epsilon\searrow 0}\frac 1\epsilon \mathbb{P}_{(C\cup \{i\})^{\rm c},+,0}(|X_{i}|\leq
\epsilon) \geq \lim_{\epsilon\searrow 0}\frac 1\epsilon \mathbb{P}_{\mathbb{Z}^2,+,0}(|X_{i}|\leq
\epsilon \,\vert\, X_j=0)\,,
\end{equation*}
where $j$ in the last probability is a site of $C$ which is nearest
neighbor to $i$ (the last inequality follows from FKG). Let us estimate this
probability. By the Markov property, we have ($u\sim v$ meaning $u$ and $v$
nearest neighbors)
\begin{align*}
\mathbb{P}_{\mathbb{Z}^2,+,0}(|X_{i}|\leq
\epsilon &\,\vert\, X_j=0)\\
&\geq \int\mathbb{P}_{\mathbb{Z}^2,+,0}(\mathrm{d} X
\,|\, X_j=0)\; 1_{\{\max_{k\sim i} X_k < h_1\}}
\mathbb{P}_{\{i\},+,X}(X_i\leq\epsilon)\\
&\geq \int\mathbb{P}_{\mathbb{Z}^2,+,0}(\mathrm{d} X\,|\, X_j=0)\; 1_{\{\max_{k\sim i} X_k < h_1\}}
\mathbb{P}_{\{i\},+,0}^{h_1}(X_i\leq\epsilon)\\
&=\mathbb{P}_{\mathbb{Z}^2,+,0}( \max_{k\sim i} X_k < h_1 \,|\,
X_j=0)\;\mathbb{P}_{\{i\},+,0}^{h_1}(X_i\leq\epsilon)\,.
\end{align*}
Now, this last probability can easily be bounded from below:
\begin{equation*}
\mathbb{P}_{\{i\},+,0}^{h_1}(X_i\leq\epsilon) \geq c_1\;\epsilon\,,
\end{equation*}
while the other probability can be bounded below by $1/2$ if $h_1$ is chosen
large enough. Indeed, since
\begin{equation*}
\mathbb{E}_{\mathbb{Z}^2,+,0} (X_k \,|\, X_j = 0) \leq c_2\,,
\end{equation*}
when $k\sim i\sim j$, Markov inequality yields
\begin{equation*}
\mathbb{P}_{\mathbb{Z}^2,+,0}(\max_{k\sim i} X_k \geq h_1 \,|\,
X_j=0) \leq 4\;c_2/h_1\,.
\end{equation*}
\end{proof}

\noindent
\begin{proof}[Proof of Theorem~\ref{thm_2ndOrder}]
The proof is trivial; by FKG inequalities,
$$
\mathbb{E}_{\mathbb{Z}^2,+,0,\upsilon} (X_i) \geq \mathbb{E}_{\mathbb{Z}^2,0,\upsilon} (X_i 1_{\{X_i\geq 0\}})
= \tfrac12\mathbb{E}_{\mathbb{Z}^2,0,\upsilon} (|X_i|) \geq C\sqrt{|\log\upsilon|}\,,
$$
the last inequality being proved in~\cite{BoVe2001}.
\end{proof}
%
%
%
%
%
%
\subsection{Proof of Theorem~\ref{thm_Ising}}\label{2dIsing}
\subsubsection{Preliminaries}
Let us first introduce our notations; note that, for reasons of convenience, the
notations we use here differ slightly from those used in the introduction.

Let $\Lambda\Subset\mathbb{Z}^2$, $h\in\mathbb{R}$, and let, for each bond
$e=\langle x,y \rangle$, $x,y\in\mathbb{Z}^2$,
$$
J(e) =
\begin{cases}
  h  &  \text{if $\min (x_2,y_2) = 0$ and $\max (x_2,y_2) = 1$;}\\
  1  &  \text{otherwise.}
\end{cases}
$$
We define the following Hamiltonian acting on configurations
$\sigma\in\{-1,1\}^{\Lambda}$,
$$
H_{\lambda,h,\Lambda} (\sigma) = -\sum_{\langle x,y \rangle \subset \Lambda}
J(\langle x,y \rangle)\,
(\sigma_x\sigma_y-1) - \lambda \sum_{x\in\Lambda} \sigma_x\,.
$$
$\lambda$ and $h$ are respectively the bulk and boundary magnetic fields.
Let $\partial\Lambda = \{x\in\Lambda\,:\, \exists y\not\in\Lambda,\, |x-y|=1\}$.
We define three
Gibbs measures in $\Lambda$; all are probability measures on
$\sigma\in\{-1,1\}^{\Lambda}$. Let $s=\pm 1$;
$$
\mu_{s,\beta,\lambda,h,\Lambda} =
\begin{cases}
(Z_{s,\beta,\lambda,h,\Lambda})^{-1}\; e^{-\beta\, H_{\lambda,h,\Lambda} (\sigma)}
   & \text{if } \sigma_x = s,\, \forall x\in\partial\Lambda\,;\\
0  &  \text{otherwise.}
\end{cases}
$$
and
$$
\mu_{\pm,\beta,\lambda,h,\Lambda} =
\begin{cases}
(Z_{\pm,\beta,\lambda,h,\Lambda})^{-1}\; e^{-\beta\, H_{\lambda,h,\Lambda}
(\sigma)}
   & \text{if } \sigma_x = \mathrm{sign}(x_2),\, \forall x\in\partial\Lambda\,;\\
0  &  \text{otherwise,}
\end{cases}
$$
where we set $\mathrm{sign}(0) = -1$.

Let $\beta>\beta_c$. We denote by $h_w(\beta)\geq 0$ the value of the boundary
field at which wetting takes place when $\lambda=0$.
Since $\beta>\beta_c$ and $h\geq h_w(\beta)$ are kept fixed, we often lighten
the notations by omitting the corresponding subscripts.

For $N$ even, consider the square box
\[
\Lambda_N=\{-N/2+1,\dots,N/2\} \times \{0,\dots,N-1\}.
\]

In the complete wetting regime $h\geq h_w(\beta)$, when $\lambda=0$ the wall free energy is given by the surface tension in the horizontal direction. In the presence of a bulk field $\lambda>0$, the latter does not make sense anymore since the $-$ phase isn't stable, however the former remains meaningful since in the vicinity of the wall the $-$ phase is stabilized by the boundary conditions. We therefore define the (finite-volume) wall free energy by the usual formula,
$$
\tau_{\mathrm{bd}}(\beta,h,\lambda,N) = - \frac1N\, \log\frac{Z_{\pm,\beta,\lambda,h,\Lambda_N}}{Z_{+,\beta,\lambda,h,\Lambda_N}}\,.
$$
The next lemma states that the effect of the bulk field $\lambda$
on the surface tension is to increase the latter by an amount of order at most
$\lambda^{2/3 + o(1)}$.
\begin{lemma}\sl
\label{lem_lowerboundPF_Ising}
For all $\beta>\beta_c$ and $h>h_w(\beta)$, there exist $C>0$, $\lambda_0>0$
and $K>0$ such that, for all $0<\lambda<\lambda_0$ and
$N>K\lambda^{-2/3}|\log\lambda|^{3}$,
$$
\tau_{\mathrm{bd}}(\beta,h,0,N)
\leq \tau_{\mathrm{bd}}(\beta,h,\lambda,N)
\leq \tau_{\mathrm{bd}}(\beta,h,0,N) + C \lambda^{2/3}|\log\lambda|^2\,.
$$
\end{lemma}

\begin{proof}[Proof of Lemma~\ref{lem_lowerboundPF_Ising}]

Let $F_\lambda(\sigma)=e^{\lambda\,\sum_{x\in\Lambda_N} \sigma_x}$.
The statement of the lemma can be rewritten in the following form,
$$
\langle F_\lambda \rangle_{+,\beta,0,h,\Lambda_N}
\geq \langle F_\lambda \rangle_{\pm,\beta,0,h,\Lambda_N}
\geq e^{- C_4\,\lambda^{2/3}\, |\log\lambda|^2\, N}
\langle F_\lambda \rangle_{+,\beta,0,h,\Lambda_N}\,.
$$
Since the first inequality is a direct consequence of FKG inequalities, we only have to prove the second one.

Let $\mathcal{H}=\lambda^{-1/3}|\log\lambda|^2$ and $\mathcal{S}_\mathcal{H} = \{x\in\Lambda_N
\,:\, x_2 \leq \mathcal{H}\}$ be the strip of width $\mathcal{H}$ along the bottom wall.

According to the discussion in Appendix~\ref{app_randomline},
$$
Z_{\pm,\beta,\lambda,h,\Lambda_N} \geq \sum_{\gamma\subset \mathcal{S}_\mathcal{H}}
w(\gamma)\,
Z_{+,\beta,\lambda,\Lambda^+(\gamma)}\, Z_{-,\beta,\lambda,h,\Lambda^-(\gamma)}\,.
$$
To simplify somewhat the notations we omit, in the rest of the proof, the
subscripts $\beta$ and $h$, and write simply $\Lambda^+$ and $\Lambda^-$.

Our first task is to remove the $\lambda$-dependence in the last equation.
Indeed our main tool, the random-line representation briefly described in
Appendix~\ref{app_randomline}, only applies in the absence of bulk field.
\begin{align*}
Z_{+,\lambda,\Lambda^+}\,Z_{-,\lambda,\Lambda^-}
&\geq Z_{+,\lambda,\Lambda^+}\,Z_{-,-\lambda,\Lambda^-}
\;e^{-2\lambda\, |\Lambda^-(\gamma)|}\\
&\geq \frac{Z_{+,\lambda,\Lambda^+}\,
Z_{-,-\lambda,\Lambda^-}} {Z_{+,0,\Lambda^+}\,
Z_{+,0,\Lambda^-}} \; Z_{+,0,\Lambda^+}\,
Z_{+,0,\Lambda^-}\;e^{-2\lambda\, \mathcal{H}\,N}\,.
\end{align*}
Writing $x\sim\gamma$ if $x$ is one the sites whose value is completely
determined by $\gamma$, we can estimate the ratio
\begin{align*}
\frac{Z_{+,\lambda,\Lambda^+}\,
Z_{-,-\lambda,\Lambda^-}} {Z_{+,0,\Lambda^+}\,
Z_{+,0,\Lambda^-}}
&= \frac{\mu_{+,\lambda,\Lambda_N}(\sigma_x = 1,\, \forall
x\sim\gamma)} {\mu_{+,0,\Lambda_N}(\sigma_x = 1,\, \forall
x\sim\gamma)}\; \frac{Z_{+,\lambda,\Lambda_N}}{Z_{+,0,\Lambda_N}}\\
&\geq  \frac{Z_{+,\lambda,\Lambda_N}}{Z_{+,0,\Lambda_N}}\,.
\end{align*}
Here, the identity follows from the $+/-$ symmetry of the model, and the last
inequality is a consequence of FKG inequalities.
Collecting all these estimates together, we get, see~\eqref{def_q},
$$
Z_{\pm,\lambda,\Lambda_N} \geq e^{- 2\lambda\, \mathcal{H}\, N}\, Z_{+,\lambda,\Lambda_N} \,
\sum_{\gamma\subset\mathcal{S}_\mathcal{H}} q_{N}(\gamma)\,,
$$
which is precisely what we were after. At this stage, it is possible to use the
tools discussed in Appendix~\ref{app_randomline} in order to control the last
sum.

Let $K$ be a sufficiently large integer.
We split the slab $\mathcal{S}_\mathcal{H}$ into $M$
disjoint rectangles, $\mathcal{R}_1,\dots,\mathcal{R}_M$,
with height $\mathcal{H}$ and basis of length $[\mathcal{H}^2/(K\log \mathcal{H})]$, except possibly for the
rightmost one which may have a shorter basis.
We also introduce the dual sites $a_k$, $k=0,\dots,M$, defined by
$$
a_k = (-(N-1)/2+k[\mathcal{H}^2/(K\log \mathcal{H})] \wedge (N-1)/2,1/2)\,.
$$
($a_{k-1}$ and $a_k$ are thus the bottom corners of the rectangle
$\mathcal{R}_k$).

\begin{figure}[t!]
\centerline{\resizebox{!}{2cm}{\input{isinglayer.pstex_t}}}
\caption{The restricted family of open contours.}
\label{fig_layer}
\end{figure}

We further restrict the summation to the set of open contours $\gamma$
satisfying the following conditions (see Fig.~\ref{fig_layer}):
\begin{itemize}
\item $\gamma = \eta_1 \amalg \eta_2 \amalg \dots \amalg \eta_M$.
\item The piece $\eta_k$, $k=1,\dots,M$, use only inner edges of the $k$th
rectangle, and connects $a_{k-1}$ and $a_k$.
\end{itemize}
By~\eqref{infvol} and~\eqref{split},
$$
q_{N}(\gamma)\geq\prod_{k=1}^M q_{\mathrm{s.i.}}(\eta_k)\,.
$$
Therefore
$$
\sum_{\gamma\subset\mathcal{S}_\mathcal{H}} q_{N}(\gamma) \geq \prod_{k=1}^M
\sum_{\eta_k\subset\mathcal{R}_k} q_{\mathrm{s.i.}}(\eta_k)\,;
$$
the last summation runs over contours satisfying the condition above. It is
controlled thanks to~\eqref{Concentration}, which implies that
$$
\sum_{\eta_k\subset\mathcal{R}_k} q_{\textrm{s.i.}}(\eta_k) \geq \mathcal{H}^{-C}\;
\sum_{\zeta:\, a_{k-1} \to a_k} q_{\textrm{s.i.}}(\zeta)\,,
$$
provided $K$ is big enough. We then combine~\eqref{exp2PF} and~\eqref{OZ32} to
obtain
\begin{align*}
\sum_{\zeta:\, a_{k-1} \to a_k} q_{\textrm{s.i.}}(\zeta)
&\geq
C\,|a_{k-1}-a_k|^{-3/2}\,e^{-\tau_\beta(\mathbf{e}_1)|a_{k-1}-a_k|}\\
&\geq
C\, \mathcal{H}^{-3}\, e^{-\tau_\beta(\mathbf{e}_1)\, |a_{k-1}-a_k|}\,,
\end{align*}
where $\tau_\beta(\mathbf{e}_1)$ is the surface tension in the horizontal
direction. Since, by~\eqref{expansion} (with $\lambda=0$) and~\eqref{def_q} ,
\eqref{exp2PF} and~\eqref{OZ32},
$$
e^{-\tau_\beta(\mathbf{e}_1)\, N} \geq
\frac{Z_{\pm,0,\Lambda_N}}{Z_{+,0,\Lambda_N}}\,,
$$
the conclusion follows.
\end{proof}
%
%
\subsubsection{Proof of the lower bound}
We turn now to the proof of~\eqref{eq_Ising_LB}.
As before, we omit reference to $\beta$ and $h$ in the notations,
write simply  $\Lambda^+$ and $\Lambda^-$, and also set $C_{\lambda,N} =
\lambda^{-1/3}|\log\lambda|^{-3}\, N$. We also denote by $\mathcal{A}$ the set
of open contours such that $|\Lambda^-| < C_{\lambda,N}$. Applying
Lemma~\ref{lem_lowerboundPF_Ising} yields
\begin{equation*}
\mu_{\pm,\lambda,\Lambda_N} (\mathcal A)
\leq
e^{C\,\lambda^{2/3}|\log\lambda|^2\, N}\; \frac{Z_{+,0,\Lambda_N}}
{Z_{\pm,0,\Lambda_N}}\; \sum_{\gamma\in\mathcal{A}} w(\gamma)\;
\frac{Z_{+,\lambda,\Lambda^+}\, Z_{-,\lambda,\Lambda^-}}
{Z_{+,\lambda,\Lambda_N}}\,.
\end{equation*}
Our first task is again to get rid of the bulk field. This is the content of the
following
%
%
\begin{lemma}\label{lem_OpenContours}
For all $\beta>\beta_c$, there exists $C_1<\infty$ and $C_2>0$ such that, for
any set $\mathcal{C}$ of open contours,
\begin{multline}
\label{nofield}
\frac{Z_{+,0,\Lambda_N}} {Z_{\pm,0,\Lambda_N}}\;
\sum_{\gamma\in\mathcal{C}} w(\gamma)\;
\frac{Z_{+,\lambda,\Lambda^+}\, Z_{-,\lambda,\Lambda^-}}
{Z_{+,\lambda,\Lambda_N}}\\
\leq e^{C_1\, \lambda\, N}\; \frac{Z_{+,0,\Lambda_N}} {Z_{\pm,0,\Lambda_N}}\;
\sum_{\gamma\in\mathcal{C}} q_{N}(\gamma)\; + \; e^{-C_2\,N}\,.
\end{multline}
\end{lemma}
\begin{proof}[Proof of Lemma~\ref{lem_OpenContours}]
By symmetry,
\begin{align*}
\frac{Z_{+,\lambda,\Lambda^+}\, Z_{-,\lambda,\Lambda^-}}
{Z_{+,\lambda,\Lambda_N}}
&= \frac{Z_{+,-\lambda,\Lambda^-}}
{Z_{+,\lambda,\Lambda^-}}\;\mu_{+,\lambda,\Lambda_N} (\sigma_x=1,\, \forall
x\sim\gamma)\nonumber\\
&\leq\mu_{+,\lambda,\Lambda_N} (\sigma_x=1,\, \forall
x\sim\gamma)\,,
\end{align*}
where as before $x\sim\gamma$ if its value is completely determined by $\gamma$,
and we used
\begin{equation*}
\frac{Z_{+,-\lambda,\Lambda^-}} {Z_{+,\lambda,\Lambda^-}} = \exp\left\{
\sum_{x\in\Lambda^-} \int_0^\lambda \left( \langle
\sigma_x \rangle_{+,-s,\Lambda^-} -  \langle
\sigma_x \rangle_{+,s,\Lambda^-}\right)\; \mathrm{d}s \right\} \leq 1\,.
\end{equation*}
Now,
\begin{multline*}
\mu_{+,\lambda,\Lambda_N} (\sigma_x=1,\, \forall
x\sim\gamma) =
\mu_{+,0,\Lambda_N} (\sigma_x=1,\, \forall
x\sim\gamma)\\
\times \exp\left\{ \sum_{x\in\Lambda_N} \int_0^\lambda \bigl( \langle \sigma_x
\,|\, \sigma_y=1,\, \forall y\sim\gamma \rangle_{+,s,\Lambda_N} -  \langle
\sigma_x \rangle_{+,s,\Lambda_N}\bigr)\; \mathrm{d}s \right\}\,.
\end{multline*}
The difference in the exponent can easily be bounded using~\eqref{BLP}:
There exists $C<\infty$ such that, uniformly in $s\geq 0$,
$$
\langle \sigma_x \,|\, \sigma_y=1,\, \forall y\sim\gamma
\rangle_{+,s,\Lambda_N} -  \langle \sigma_x \rangle_{+,s,\Lambda_N}
\leq C\; \sum_{z\sim\gamma} e^{-|x-z|/C}\,,
$$
and therefore
$$
\sum_{x\in\Lambda} \int_0^\lambda \bigl( \langle \sigma_x \,|\, \sigma_y=1,\,
\forall y\sim\gamma \rangle_{+,s,\Lambda_N} -  \langle \sigma_x
\rangle_{+,s,\Lambda_N}\bigr)\, \mathrm{d}s \leq C'\,\lambda\, \#\{z\,:\,
z\sim\gamma\}\,.
$$
We are almost done. Observe now that $\#\{z\,:\, z\sim\gamma\} < 4|\gamma|$. Let $K$ be a large number, which will be
chosen below. In the case $|\gamma| < K\,N$, we get
$$
\mu_{+,\lambda,\Lambda_N} (\sigma_x=1,\, \forall
x\sim\gamma) \leq e^{4C'K\, \lambda\,N}\;
\mu_{+,0,\Lambda_N} (\sigma_x=1,\, \forall
x\sim\gamma)\,.
$$
This yields the first term in~\eqref{nofield}. The second term takes care of
the remaining open contours. Indeed,
\begin{align*}
\sum_{\substack{\gamma:\\|\gamma|\geq KN}} q_{N}(\gamma)\;
e^{CC'\,\lambda\,|\gamma|}
&\leq \sum_{k\geq 0} e^{CC'\,\lambda\, (K+k+1)\, N}\; \sum_{l=(K+k)\,
N}^{(K+k+1)\, N -1} \sum_{\substack{\gamma:\\|\gamma|= l}} q_{N}(\gamma)\\
&\leq \sum_{k\geq 0} e^{-C''(K+k)\,N}\\
&\leq e^{-C'''K\,N}\,,
\end{align*}
where we used~\eqref{longunlikely}. Consequently, since
$$
\frac{Z_{+,0,\Lambda_N}} {Z_{\pm,0,\Lambda_N}} \leq
C\, N^{3/2}\, e^{\tau_\beta(\mathbf{e}_1)\, N}\,,
$$
we see that
$$
\frac{Z_{+,0,\Lambda_N}} {Z_{\pm,0,\Lambda_N}}\;
\sum_{\gamma:\, |\gamma| \geq K\,N} w(\gamma)\;
\frac{Z_{+,\lambda,\Lambda^+}\, Z_{-,\lambda,\Lambda^-}}
{Z_{+,\lambda,\Lambda_N}}
\leq e^{-C\, N}\,,
$$
provided $K$ is chosen large enough.
\end{proof}

Applying the lemma,  we have
\begin{align*}
\mu_{\pm,\lambda,\Lambda_N} (\mathcal A)
&\leq e^{C\,\lambda^{2/3}|\log\lambda|^2\, N}\;
\frac{Z_{+,0,\Lambda_N}} {Z_{\pm,0,\Lambda_N}}\;
\sum_{\substack{\gamma:\\|\Lambda^-| < C_{\lambda,N}}} q_{N}(\gamma) +
e^{-a\,N}\\
&= e^{C\,\lambda^{2/3}|\log\lambda|^2\, N}\;
\mu_{\pm,0,\Lambda_N} (|\Lambda^-| < C_{\lambda,N}) +
e^{-a\,N}\,,
\end{align*}
for some $a>0$, and so we are left with estimating the probability that
$|\Lambda^-| < C_{\lambda,N}$ in the absence of bulk magnetic field, a purely
entropic problem.

\begin{figure}[t!]
\centerline{\resizebox{!}{3.5cm}{\input{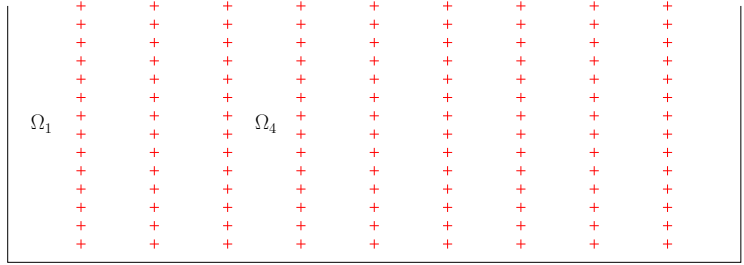}}}
\caption{The box $\Omega$ is obtained by splitting $\Lambda_N$ into strips of
width $[A\lambda^{-2/3}|\log\lambda|^{-2}]$ by forcing all the spins on the
corresponding vertical half-lines to take the value $+1$.}
\label{fig_Omega}
\end{figure}

The next step is to observe that the event $\{|\Lambda^-| <
C_{\lambda,N}\}$ is increasing. So, by FKG inequalities,
$$
\mu_{\pm,0,\Lambda_N} (|\Lambda^-| < C_{\lambda,N}) \leq
\mu_{\pm,0,\Omega} (|\Lambda^-| < C_{\lambda,N})\,
$$
where $\Omega$ is the box (see Fig.~\ref{fig_Omega})
$$
\Omega = \Lambda_N \setminus \{z:\, z_1 = -N/2 + 1 +
k[A\lambda^{-2/3}|\log\lambda|^{-2}],\,
k\in\mathbb{N}\}\,,
$$
$A$ being a small number to be chosen later.
Let us denote by $\Lambda^-_k$ the component of $\Lambda^-$ contained in the
$k$th slice $\Omega_k$ of $\Omega$. Let $K_{\lambda} =
\lambda^{-1}|\log\lambda|^{-9/2}$. We are going to show that
\begin{equation}
\mu_{\pm,0,\Omega_1}(|\Lambda^-_1| < K_{\lambda} ) < 7/8\,.
\label{7_8}
\end{equation}
From this and a standard large deviations estimate, we get that
$$
\mu_{\pm,0,\Omega} \left( \# \{ k:\, |\Lambda^-_k| < K_{\lambda} \} \geq
\frac{15}{16}\, \frac{N}{[A\lambda^{-2/3}|\log\lambda|^{-2}]}\right) \leq
e^{-CA^{-1}\, \lambda^{2/3}\, |\log\lambda|^2\, N}\,.
$$
(Notice that the events $\{ |\Lambda^-_k| < K_{\lambda} \}$ are
independent under $\mu_{\pm,0,\Omega}$). The proposition easily follows from
this, provided $\lambda$ is chosen small enough (for a fixed small $A$). Indeed, on the complementary event,
i.e.
$$
\# \{ k:\, |\Lambda^-_k| >
K_{\lambda} \} > \frac{1}{16}\,
\frac{N}{[A\lambda^{-2/3}|\log\lambda|^{-2}]}\,,
$$
we have
$$
|\Lambda^-| > \tfrac1{16} N/[A\lambda^{-2/3}|\log\lambda|^{-2}]
\,K_{\lambda} > C_{\lambda,N}\,,
$$
for small $\lambda$. Let us now prove~\eqref{7_8}.

\begin{figure}[t!]
\centerline{\resizebox{!}{5cm}{\input{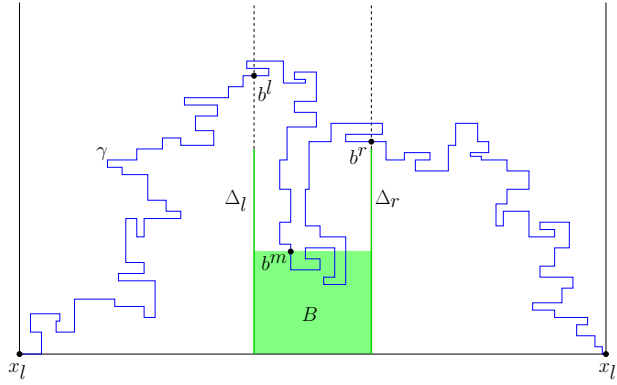}}}
\caption{The settings for the proof of~\eqref{7_8}.}
\label{fig_minvol}
\end{figure}

Denote by $x^l$ and $x^r$ the dual sites which are the bottom left and bottom
right corners of $\Omega_1$ (so these are the endpoints of the open path in
$\Omega_1$, see Fig.~\ref{fig_minvol}). Let $a_l = x^l_1 + [(\tfrac12 -
A^2|\log\lambda|^{-1})A\lambda^{-2/3}|\log\lambda|^{-2}]$ and $a_r = x_r -
[(\tfrac12 - A^2 |\log\lambda|^{-1})A\lambda^{-2/3}|\log\lambda|^{-2}]$. We
introduce the following sets of dual sites:
\begin{align*}
\Delta_l &= \{ x \in \mathbb{Z}^{2,*}: \, x_1 = a_l,\, \tfrac12 \leq x_2 < A
\lambda^{-1/3} |\log\lambda|^{-1} \}\\
\Delta_r &= \{ x \in \mathbb{Z}^{2,*}: \, x_1 = a_r,\, \tfrac12 \leq x_2 < A
\lambda^{-1/3} |\log\lambda|^{-1} \}\\
\Delta &= \Delta_l\cup\Delta_r\\
B &= \{ x \in \mathbb{Z}^{2,*}: \, a_l < x_1 < a_r,\, \tfrac12 \leq x_2 < \tfrac12 A
\lambda^{-1/3} |\log\lambda|^{-1} \}\,.
\end{align*}
We are going to show that
\begin{equation}
\sum_{\substack{\gamma:x^l\to x^r\\\gamma\cap B = \emptyset}} q_{\Omega_1}
(\gamma) \geq \tfrac18\, \langle \sigma_{x^l} \sigma_{x^r}
\rangle_{\Omega_1}\,,
\label{GoesHigh}
\end{equation}
which readily implies~\eqref{7_8} since $|B| > K_{\lambda}$ (when $\lambda$ is
small). To prove~\eqref{GoesHigh}, we first show that the open contour
typically does not hit $\Delta$. This is equivalent to saying that
$\Lambda_1^-$ contains the set $\Delta$. But this is a decreasing event.
Therefore, if we introduce the new box
$$
\widetilde\Omega_1 = \{x=(x_1,x_2):\, (x_1,|x_2|) \in \Omega_1\}\,,
$$
and the new set
$$
\widetilde\Delta = \{x=(x_1,x_2):\, (x_1,|x_2|) \in \Delta\}\,,
$$
then by FKG inequalities,
$$
\mu_{\pm,0,\Omega_1}(\Lambda^-_1 \supset \Delta) \geq
\mu_{\pm,0,\widetilde\Omega_1}(\Lambda^-_1 \supset \Delta)
\geq \mu_{\pm,0,\widetilde\Omega_1}(\Lambda^-_1 \supset \widetilde\Delta)\,.
$$
The latter probability is easily bounded. Indeed, by symmetry it is larger than
$$
\tfrac12\,\langle \sigma_{x^l} \sigma_{x^r} \rangle_{\widetilde\Omega_1}^{-1}\;
\sum_{\substack{\gamma:x^l\to x^r\\\gamma\cap \widetilde\Delta =
\emptyset}} q_{\widetilde\Omega_1} (\gamma),\,
$$
and the latter expression is larger than $1/4$ since
\begin{align*}
\sum_{\substack{\gamma:x^l\to x^r\\\gamma\cap \widetilde\Delta \neq
\emptyset}} q_{\widetilde\Omega_1} (\gamma)
&\leq \sum_{u\in\widetilde\Delta} \langle \sigma_{x^l}
\sigma_{u} \rangle_{\widetilde\Omega_1}\, \langle \sigma_{u}
\sigma_{x^r} \rangle_{\widetilde\Omega_1}\\
&\leq C\, \sum_{u\in\widetilde\Delta}
\frac{e^{-\tau(u-x^l)}}{|u-x_l|^{1/2}} \frac{e^{-\tau(x^r-u)}}{|x_r-u|^{1/2}}\\
&\leq 3C\, \frac{e^{-\tau(x^r-x^l)}}{|x_r-x_l|^{1/2}} \;
\frac{|\widetilde\Delta|}{|x_r-x_l|^{1/2}} \\
&\leq \tfrac12\, \langle \sigma_{x^l} \sigma_{x^r}
\rangle_{\widetilde\Omega_1}\,,
\end{align*}
provided $A$ is chosen small enough. We have thus shown that
\begin{equation}
\mu_{\pm,0,\Omega_1}(\Lambda^-_1 \supset \Delta) \geq 1/4\,.
\label{typicallyhigh}
\end{equation}

We now restrict our attention to open paths $\gamma$ in $\Omega_1$ not
intersecting $\Delta$, but intersecting $B$. To each such paths we
associate three dual sites: the first visit to the set $\{x : \, x_1 = a_l,\,
x\not\in\Delta_l\}$, the last visit to
$\{x : \, x_1 = a_r,\, x\not\in\Delta_r\}$, and the first visit to $B$. We
denote these dual sites by $b^l$, $b^r$ and $b^m$ respectively.
Using~\eqref{skel}, we then have
\begin{align*}
\sum_{\substack{\gamma:x^l\to x^r\\
                \gamma\cap \Delta = \emptyset\\
		\gamma\cap B \neq \emptyset
	       }
     } q_{\Omega_1} (\gamma)
&\leq \sum_{b^l,b^m,b^r}
\langle \sigma_{x^l} \sigma_{b^l} \rangle_{\Omega_1}\;
\langle \sigma_{b^l} \sigma_{b^m} \rangle_{\Omega_1}\;
\langle \sigma_{b^m} \sigma_{b^r} \rangle_{\Omega_1}\;
\langle \sigma_{b^r} \sigma_{x^r} \rangle_{\Omega_1}
\end{align*}
We first observe that
$$
\sum_{b^l:\, b^l_2 > \lambda^{-1}} \langle \sigma_{x^l} \sigma_{b^l}
\rangle_{\Omega_1} \leq e^{-C\, \lambda^{-1}} \ll \langle \sigma_{x^l}
\sigma_{x^r} \rangle_{\Omega_1}\,,
$$
and similarly for $b^r$, so that we can restrict the sum over these points to
those with second component smaller than $\lambda^{-1}$. We now use the sharp triangle inequality~\eqref{STI} to get
\begin{align*}
\tau(b^m-b^l) + \tau(b^r-b^m)
&\geq \tau(b^r-b^l) + \kappa\, (|b^m-b^l|+|b^r-b^m|-|b^r-b^l|)\\
&\geq \tau(b^r-b^l) + CA^{-1}\, |\log\lambda|\,,
\end{align*}
and therefore, using~\eqref{OZ},
\begin{align*}
\langle \sigma_{b^l} \sigma_{b^m} \rangle_{\Omega_1}\;
\langle \sigma_{b^m} \sigma_{b^r} \rangle_{\Omega_1}
&\leq
\langle \sigma_{b^l} \sigma_{b^m} \rangle\;
\langle \sigma_{b^m} \sigma_{b^r} \rangle\\
&\leq
e^{-C'A^{-1}\, |\log\lambda|}\;
\langle \sigma_{b^l} \sigma_{b^r} \rangle\, |b^r-b^l|\\
&\leq
e^{-C''A^{-1}\, |\log\lambda|}\;
\langle \sigma_{b^l} \sigma_{b^r} \rangle_{\Omega_1}\,.
\end{align*}
Using this and
$$
\langle \sigma_{x^l} \sigma_{b^l} \rangle_{\Omega_1}\;
\langle \sigma_{b^l} \sigma_{b^r} \rangle_{\Omega_1}\;
\langle \sigma_{b^r} \sigma_{x^r} \rangle_{\Omega_1}
\leq
\langle \sigma_{x^l} \sigma_{x^r} \rangle_{\Omega_1}\,,
$$
which follows from GKS inequalities, we finally get
$$
\sum_{\substack{\gamma:x^l\to x^r\\
                \gamma\cap \Delta = \emptyset\\
		\gamma\cap B \neq \emptyset
	       }
     } q_{\Omega_1} (\gamma)
\leq e^{-C'''A^{-1}\, |\log\lambda|}\; \langle \sigma_{x^l} \sigma_{x^r}
\rangle_{\Omega_1} \leq \tfrac18\, \langle \sigma_{x^l} \sigma_{x^r}
\rangle_{\Omega_1}\,.
$$
We deduce from the latter bound and~\eqref{typicallyhigh} that
$$
\sum_{\substack{\gamma:x^l\to x^r\\
		\gamma\cap B \neq \emptyset
	       }
     } q_{\Omega_1} (\gamma)
\leq (\tfrac18 + \tfrac34)\, \langle \sigma_{x^l} \sigma_{x^r}
\rangle_{\Omega_1}\,,
$$
and~\eqref{GoesHigh} is proved.

\subsubsection{Proof of the upper bound}
Let us write $V=K|\log\lambda|^2\lambda^{-1/3}N$, where $K$ will be chosen
sufficiently large later on. We want to bound from above the following
probability: $$
\mu_{\pm,\lambda,\Lambda_N} \left( |C^-(\,\cdot\,)| \geq V \right)\,,
$$
where $C^-(\sigma)$ is the set of all sites in $\Lambda_N$ that are
connected to the bottom wall by a path of $-$ spins in $\sigma$. Our aim is to show that this probability goes to zero with $N$, provided $K$ has been chosen large enough.

We start with the following obvious upper bound:
\begin{equation*}
\frac{Z_{\pm,\lambda,\Lambda_N}\left( |C^-(\,\cdot\,)| \geq V \right)}
     {Z_{\pm,\lambda,\Lambda_N}}
\leq
\frac{Z_{\pm,\lambda,\Lambda_N}\left( |C^-(\,\cdot\,)| \geq V \right)}
     {Z_{\pm,0,\Lambda_N}\left( |C^-(\,\cdot\,)| \geq V \right)}
\,
\frac{Z_{\pm,0,\Lambda_N}}
     {Z_{\pm,\lambda,\Lambda_N}}
\end{equation*}
By Lemma~\ref{lem_lowerboundPF_Ising},
$$
\frac{Z_{\pm,0,\Lambda_N}}
     {Z_{\pm,\lambda,\Lambda_N}}
\leq e^{C|\log\lambda|^2\lambda^{2/3}N}
\frac{Z_{+,0,\Lambda_N}}
     {Z_{+,\lambda,\Lambda_N}}\,.
$$
Now observe that we can write
$$
\frac{Z_{+,0,\Lambda_N}}
     {Z_{+,\lambda,\Lambda_N}}
=
\exp\left\{ - \int_0^\lambda \mathrm{d}\lambda' \sum_{x\in\Lambda_N} \langle \sigma_x
\rangle_{+,\lambda',\Lambda_N}\right\}\,,
$$
and, similarly,
$$
\frac{Z_{\pm,\lambda,\Lambda_N}\left( |C^-(\,\cdot\,)| \geq V \right)}
     {Z_{\pm,0,\Lambda_N}\left( |C^-(\,\cdot\,)| \geq V \right)}
=
\exp\left\{ \int_0^\lambda \mathrm{d}\lambda' \sum_{x\in\Lambda_N} \langle \sigma_x
\,|\, |C^-(\,\cdot\,)| \geq V \rangle_{\pm,\lambda',\Lambda_N}\right\}\,.
$$
We estimate the product of these two expressions using a coupling argument.
Since, by FKG inequalities and the fact that $\{ |C^-(\,\cdot\,)| \geq V \}$ is
non-increasing,
$$
\mu_{+,\lambda',\Lambda_N} \succeq
\mu_{\pm,\lambda',\Lambda_N}(\,\cdot\,|\, |C^-(\,\cdot\,)| \geq V)\,,
$$
there exists a coupling $\nu_{\lambda'}$ of $\mu_{+,0,\Lambda_N}$ and
$\mu_{\pm,\lambda',\Lambda_N}(\,\cdot\,|\, |C^-(\,\cdot\,)| \geq V)$ such that
$$
\nu(\{\sigma\geq \sigma'\}) = 1\,.
$$
In fact, it is well-known that one can construct this coupling explicitly. We
quickly sketch this, because we'll need some non-degeneracy property below.
First, we order all the sites in $\Lambda_N$, say $i_1,\ldots,i_{|\Lambda_N|}$.
Then we introduce a family $(X_{i_k})_{k=1}^{|\Lambda_N|}$ of iid random
variables, with uniform distribution over $[0,1]$. Now we set $\sigma_{i_1}=1$ (resp.
$\sigma'_{i_1}=1$) if $\mu_{+,\lambda',\Lambda_N}(\sigma_{i_1}=1)\geq X_{i_1}$ (resp. if
$\mu_{\pm,\lambda',\Lambda_N}(\sigma_{i_1}=1 \,|\, |C^-(\,\cdot\,)| \geq V) \geq
X_{i_1}$); otherwise it is set to $-1$. Suppose we have already constructed the
first $k-1$ spins. Then we set $\sigma_{i_k}=1$ if
$\mu_{+,\lambda',\Lambda_N}(\sigma_{i_k}=1 \,|\, \sigma_{i_l},\, l<k) \geq X_{i_k}$,
and similarly for $\sigma'_{i_k}$. Just observe that, by FKG inequalities,
$$
\mu_{+,\lambda',\Lambda_N}(\sigma_{i_k}=1 \,|\, \sigma_{i_l},\, l<k) \geq
\mu_{\pm,\lambda',\Lambda_N}(\sigma_{i_k}=1 \,|\, |C^-(\,\cdot\,)| \geq V,\, \sigma_{i_l},\, l<k)\,.
$$

Using this notation, we are left with estimating 
$
\sum_{x\in\Lambda_N} \left( \nu_{\lambda'}(\sigma'_x - \sigma_x) \right)\,.
$
Since $\nu_{\lambda'}(\sigma\geq \sigma')=1$, we can write
\begin{align*}
\nu_{\lambda'}\bigl( \sum_{x\in\Lambda_N} (\sigma'_x - \sigma_x) \bigr)
&\leq
\nu_{\lambda'}\bigl( \sum_{x\in C^-(\sigma')} (\sigma'_x - \sigma_x) \bigr)\,.
\end{align*}
Now observe that $\sum_{x\in C^-(\sigma')} \sigma'_x = -|C^-(\sigma')|$.
Moreover, from the above construction, we see that the coupling satisfies
$$
\inf_{x\in\Lambda_N}\nu_{\lambda'}\bigl(\sigma_x \,\bigm|\, \sigma_y,\, \forall y\in\Lambda_N \setminus\{x\},
\sigma'_z,\, \forall z\in\Lambda_N \bigr) > -1-\epsilon_\beta\,,
$$
with $\epsilon_\beta$ depending only on the parameter $\beta$, for all $\lambda$
small enough. Therefore, putting everything together, we obtain
$$
\mu_{\pm,\lambda,\Lambda_N} \left( |C^-(\,\cdot\,)| \geq V \right)
\leq
e^{-( K\epsilon_\beta - c ) |\log\lambda|^2\lambda^{2/3} N}\,,
$$
and the conclusion follows once $K$ is taken large enough.

%
%
%
%
%
%
\begin{appendix}
\section{Some tools}\label{app_randomline}
In this section we collect several results from~\cite{PfVe1997a,Ve1997,PfVe1999a} which we
are using in the paper. This is not supposed to be an introduction to the
random-line representation, and we refer the reader to the latter works for
detailed explanations.

\begin{figure}[t!]
\centerline{\includegraphics[height=35mm]{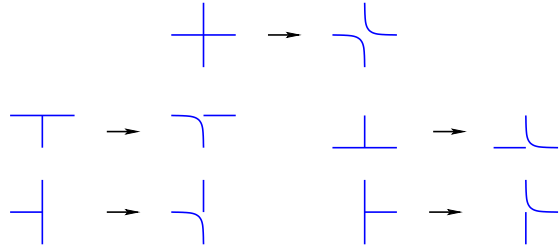}}
\caption{Deformation rules.}
\label{fig_rules}
\end{figure}

Let $\mathcal{E}^*$ be the set of all dual bonds of the infinite lattice.
One of the basic objects in this paper is the partition function
$Z_{\pm,\beta,\lambda,h,\Lambda_N}$.
To any configuration $\sigma$ with
non-zero weight under $\mu_{\pm,\beta,\lambda,h,\Lambda_N}$, we can associate a
subset $n(\sigma)$ of
$$
\mathcal{E}^*_N = \{ e^*\in\mathcal{E}^*\,:\, e^* \text{ dual to }
e\subset\Lambda_N,\, e\not\subset\partial\Lambda_N \}
$$
composed of all dual edges $e^*$ such that $e=\langle x,y
\rangle\not\subset\partial\Lambda_N$ with $\sigma_x\neq\sigma_y$. We then split
this subset into a collection of contours, by applying the rules given in
Fig.~\ref{fig_rules} each time more than two bonds are incident on a given
vertex; contours are the families of edges remaining connected after these
operations. Among all the contours, there is a unique open contour connecting
the dual sites $(-(N-1)/2,1/2)$ and $((N+1)/2,1/2)$; we call it $\gamma$.

We want to expand the partition functions in terms of the open contour $\gamma$.
Notice that once $\gamma$ is fixed, some
spins in the box $\Lambda_N$ become
frozen (i.e. in any configuration with positive weight having $\gamma$ as open
contour, these spins take the same values). We denote by $\Lambda^+(\gamma)$
the set of all non-frozen spins located
inside a component surrounded by $+$ spins. Similarly, we write
$\Lambda^-(\gamma)$ for the set of all non-frozen spins
located inside a component surrounded by $-$ spins.

We then have
\begin{equation}
\label{expansion}
Z_{\pm,\beta,\lambda,h,\Lambda_N} = \sum_{\gamma}w(\gamma)\;
Z_{+,\beta,\lambda,\Lambda^+(\gamma)}\,Z_{-,\beta,\lambda,h,\Lambda^-(\gamma)}\,,
\end{equation}
where $Z_{+,\beta,\lambda,\Lambda^+(\gamma)}$ and
$Z_{-,\beta,\lambda,h,\Lambda^-(\gamma)}$ are the two partition functions
associated to the restriction of our system to $\Lambda^+(\gamma)$ and
$\Lambda^-(\gamma)$, and $w(\gamma) = e^{-2\beta|\gamma|}$, where $|\gamma| =
\sum_{e^*\subset\gamma} J(e)$.

We introduce the weight (notice that it is computed at zero bulk field)
\begin{equation}
\label{def_q}
q_{\beta,h,N} (\gamma) = w(\gamma)\;
\frac{Z_{+,\beta,0,\Lambda^+(\gamma)}\,
Z_{-,\beta,0,h,\Lambda^-(\gamma)}} {Z_{+,\beta,0,h,\Lambda_N}}\,.
\end{equation}

Let $\beta^*$, $h^*$ and $J^*(e^*)$ be such that
$
\tanh(\beta^*) = \exp(-2\beta)$, $\tanh(\beta^*h^*) = \exp(-2\beta h)$ and
$\tanh(\beta^*J^*(e^*)) = \exp(-2\beta J(e))$.

The Gibbs measure in $\mathcal{E}^*_N$ with free b.c. is defined by
$$
\mu_{\beta^*,h^*,N}(\sigma) = (Z_{\beta^*,h^*,N})^{-1}\prod_{e^*=\langle x,y
\rangle \in\mathcal{E}^*_N} e^{-\beta^* J^*(e^*)\, \sigma_x\sigma_y}\,.
$$
A basic duality argument, see e.g.~Lemma 6.2 in~\cite{PfVe1999a}, shows that
\begin{equation}
\label{exp2PF}
\sum_{\gamma} q_{\beta,h,N} (\gamma) = \langle \sigma_{(-(N-1)/2,1/2)}
\sigma_{((N+1)/2,1/2)} \rangle_{\beta^*,h^*,N}\,.
\end{equation}
One has in fact a deeper relation between high and low temperature models.
To each subset of $\mathcal{E}^*_N$ which remains connected after applying the
rules of Fig.~\ref{fig_rules}, we can associate a weight $q^*_{\beta^*,h^*,N}$
(see~(6.8) of~\cite{PfVe1999a} for a definition) such that
$$
q^*_{\beta^*,h^*,N}(\gamma) = q_{\beta,h,N} (\gamma)\,.
$$
However, $q^*_{\beta^*,h^*,N}$ being defined for much more general subsets of
$\mathcal{E}^*_N$ than $q_{\beta,h,N}$, it is a much more useful quantity.
We'll state its properties for arbitrary $E\Subset\mathcal{E}^*$, and
will denote the corresponding weight by $q^*_{\beta^*,h^*,E}$.
We
list now the properties of these weights that we are using in the proof. Let us
denote by $\mathcal{E}^*_{\rm s.i.} = \bigcup_{N\to\infty} \mathcal{E}^*_N$.
\begin{itemize}
\item (\cite{PfVe1999a}, Lemma~6.3) The limiting weights
$q^*_{\beta^*,h^*} = \lim_{E\uparrow\mathcal{E}^*} q^*_{\beta^*,h^*,E}$
and
$q^*_{\beta^*,h^*,{\rm s.i.}} = \lim_{E\uparrow\mathcal{E}^*_{\rm s.i.}}
q^*_{\beta^*,h^*,E}$
are well-defined quantities.
Moreover, for any $\gamma\subset E$,
\begin{equation}\label{split}
q^*_{\beta^*,h^*,E}(\gamma) \geq q^*_{\beta^*,h^*,{\rm s.i.}}(\gamma) \geq
q^*_{\beta^*,h^*}(\gamma)\,.
\end{equation}
\item (\cite{PfVe1999a}, Lemma~6.4) If $\gamma = \gamma_1 \amalg \gamma_2$ ($\amalg$
is a concatenation operation, see~\cite{PfVe1999a} for a definition), then
\begin{equation}\label{infvol}
q^*_{\beta^*,h^*,E}(\gamma) \geq q^*_{\beta^*,h^*,E}(\gamma_1)\,
q^*_{\beta^*,h^*,E}(\gamma_2)\,.
\end{equation}
\item (\cite{PfVe1999a}, Lemma~6.5) Let $t_1,\dots,t_n$ be distinct dual sites in $E$.
Then
\begin{equation}\label{skel}
\sum_{\substack{\gamma:\,t_1\to\dots\to t_n}} q^*_{\beta^*,h^*,E} (\gamma) \leq
\prod_{k=1}^{n-1} \sum_{\substack{\gamma:\,t_k\to\dots\to t_{k+1}}}
q^*_{\beta^*,h^*,E} (\gamma)\,.
\end{equation}
\item (\cite{Ve1997}, Lemma 4.4.6, and \cite{PfVe1999a}, Lemma 6.10) Let $\mathcal{R}$
be a rectangular subset of $\mathcal{E}^*$ having length $R^2/(K\log R)$ and
height $R$, with basis contained inside $\{e^*=\langle x,y \rangle \in
\mathcal{E}^*\,:\, x_2=y=2=1/2\}$. We denote by $u$ and $v$ the dual sites at
the bottom left and bottom right corners of $\mathcal{R}$. Then, for $K$ large
enough, there exists $C<\infty$ such that
\begin{equation}\label{Concentration}
\sum_{\substack{\gamma:\, u \to v\\\gamma\subset\mathcal{R}}}
q^*_{\beta^*,h^*}(\gamma) \geq R^{-C}\, \sum_{\gamma:\, u \to v}
q^*_{\beta^*,h^*}(\gamma)\,.
\end{equation}
\item (\cite{PfVe1999a}, (6.9) and Lemma~6.9) For any vertices $x,y\in E$,
\begin{equation}
\label{2ptfctN}
\sum_{\gamma:\, x \to y} q^*_{\beta^*,h^*,E}(\gamma) = \langle \sigma_x
\sigma_y \rangle_{\beta^*,h^*,E}\,.
\end{equation}
Moreover
\begin{align}
\label{2ptfctsi}
\sum_{\gamma:\, x \to y} q^*_{\beta^*,h^*,{\rm s.i.}}(\gamma) &= \langle
\sigma_x
\sigma_y \rangle_{\beta^*,h^*,{\rm s.i.}}\,,\\
\intertext{and}
\sum_{\gamma:\, x \to y} q^*_{\beta^*,h^*}(\gamma) &= \langle \sigma_x
\sigma_y \rangle_{\beta^*,h^*}\,,\\
\label{2ptfct}
\end{align}
where $\langle \,\cdot\, \rangle_{\beta^*,h^*}$ and $\langle \,\cdot\,
\rangle_{\beta^*,h^*,{\rm s.i.}}$ denote expectation w.r.t. the
infinite and semi-infinite volume Gibbs measures.
\item There exists $c>0$ and $l_0<\infty$ such that, for all $l$,
\begin{equation}\label{longunlikely}
\sum_{\substack{\gamma:\\|\gamma|\geq l}} q_{\beta,h,N}(\gamma) \leq
e^{-c(l-l_0)}\,.
\end{equation}
\end{itemize}

\begin{proof}[Proof of~\eqref{longunlikely}]
The proof is similar to that of Lemma~5.6 in~\cite{PfVe1997a}; we only sketch it.
Let $R$ be some big positive number. We make a coarse-graining of the path
$\gamma$ on the scale $R$, i.e. we set $t_1$ to be the left endpoint of
$\gamma$, $t_k$ to be the first point of $\gamma$ after $t_{k-1}$ which is
outside the square of sidelength $R$ centered on $t_{k-1}$, and the procedure
stops when one reaches the other endpoint of $\gamma$. One then first sum over
these sites $t_0,t_1,\dots,t_L$ (observing that, given $t_k$ there are at most
$C R$ choices for $t_{k+1}$, and that the total length of the piece of $\gamma$
between two consecutive points is at most $R^2$, so that $L\geq l/R^2$) and
uses~\eqref{split}, \eqref{2ptfctN}, and the upper bound in~\eqref{OZ} below.
The conclusion follows since $\tau_\beta$ is uniformly strictly positive when
$\beta>\beta_c$.
\end{proof}

We also use the following results about the asymptotic behavior of the
boundary 2-point function $\langle \sigma_x \sigma_y \rangle_{\beta^*,h^*,{\rm
s.i.}}$: Suppose that $h>h_w(\beta)$; then there exist
constants $K_1$ and $K_2$ such that, for any $x,y$ with $x_2=y_2=1/2$,
\begin{equation}
\label{OZ32}
K_1\, |x-y|^{-3/2} e^{-\tau_\beta(\mathbf{e}_1)\, |x-y|}
\geq \langle \sigma_x \sigma_y \rangle_{\beta^*,h^*,{\rm s.i.}}
\geq K_2\, |x-y|^{-3/2} e^{-\tau_\beta(\mathbf{e}_1)\, |x-y|}\,,
\end{equation}
where $\tau_\beta(\mathbf{e}_1)$ is the surface tension in the horizontal
direction.

We also need the corresponding result for the bulk 2-point. There exists $K<\infty$ such that, for any $h>h_w(\beta)$, and any pair of dual sites $x,y$ in $E$ such that the set
$$
\mathcal{S}_K(x,y) = \{ u\in \mathbb{Z}^{2,*} \,:\, \|x-u\|_2 + \|y-u\|_2 \leq \|x-y\|_2 + K\log\|x-y\|_2 \}\,,
$$
satisfies $\mathcal{S}_K(x,y) \subset \{ x \in E \,:\, x_2 \geq 3/2\}$, the following holds: There exist
constants $K_1$ and $K_2$ such that (\cite{Ve1997}, Proposition 4.6.1.)
\begin{equation}
\label{OZ}
K_1\, |x-y|^{-1/2} e^{-\tau_\beta(\frac{x-y}{|x-y|})\, |x-y|}
\geq
\langle \sigma_x \sigma_y \rangle_{\beta^*,h^*}
\geq K_2\, |x-y|^{-1/2} e^{-\tau_\beta(\mathbf{e}_1)\, |x-y|}\,.
\end{equation}

An important property of surface tension is that it satisfies the following {\sl
sharp triangle inequality} (\cite{PfVe1999a}, Theorem~2.1): There exists $\kappa>0$
such that, for any $x,y$,
\begin{equation}
\label{STI}
\tau_\beta(x) + \tau_\beta(y) \geq \tau_\beta(x+y) + \kappa (|x|+|y|-|x+y|)\,.
\end{equation}

Finally, we also need the following estimate on the relaxation of correlation
functions (\cite{BrLePf1981}): Suppose that $h$ and $\lambda$ are
non-negative; then there exists $C(\beta)>0$ and $K<\infty$ such that, for any
$\Lambda_1,\Lambda_2\subset\mathbb{Z}^2$ and $A\subset \Lambda_1\cup\Lambda_2$,
\begin{equation}
\label{BLP}
\left\vert
   \langle \sigma_A \rangle_{+,\beta,\lambda,h,\Lambda_1}
 - \langle \sigma_A \rangle_{+,\beta,\lambda,h,\Lambda_2}
\right\vert
\leq K\sum_{\substack{t\in A\\t'\in\Lambda_1\triangle\Lambda_2}} e^{- C(\beta)\,
|t'-t|}\,,
\end{equation}
where $\sigma_A = \prod_{x\in A}\sigma_x$.

\end{appendix}

\bibliographystyle{plain} 
\bibliography{V03}

\end{document}